\documentclass[preprint,11pt]{article}

\usepackage{appendix}
\usepackage{longtable}
\usepackage{booktabs}
\usepackage{mathrsfs}
\usepackage{graphicx}
\usepackage{epstopdf}
\usepackage{amsmath,amsthm,amsfonts,color,xcolor,dsfont}
\usepackage{pdfsync,float,hyperref}
\usepackage{amssymb}
\usepackage{verbatim}
\usepackage{color}
\usepackage{float}
\usepackage{subfigure}

\usepackage{algorithm}
\usepackage{algpseudocode}
\usepackage{amsmath}

\numberwithin{equation}{section}
\hypersetup{
    colorlinks=true,
    linkcolor=blue,
    filecolor=magenta,
    urlcolor=magenta,
  }

\newtheorem{remark}{Remark}[section]
\newtheorem{lemma}{Lemma}[section]
\newtheorem{theorem}{Theorem}[section]
\newtheorem{definition}{Definition}[section]

\newtheorem{proposition}{Proposition}[section]

\newcommand{\noteLi}[1]{{\color{red}#1}}

\newcommand{\cstab}[1]{{\mathrm{C}_{\mathrm{stab}}(#1)}}
\begin{document}
\title{An $hp$-Adaptive Sampling Algorithm on Dispersion Relation Reconstruction for 2D Photonic Crystals}
\date{}
\author{Yueqi Wang\thanks{Department of Mathematics, The University of Hong Kong, Pokfulam Road, Hong Kong. Email: {\tt{u3007895@connect.hku.hk}} YW acknowledges support from the Research Grants Council (RGC) of Hong Kong via the Hong Kong PhD Fellowship Scheme (HKPFS).} \quad and \quad Guanglian Li\thanks{Department of Mathematics, The University of Hong Kong, Pokfulam Road, Hong Kong. Email: {\tt{lotusli@maths.hku.hk}} GL acknowledges the support from Newton International Fellowships Alumni following-on funding awarded by The Royal Society, Young Scientists fund (Project number: 12101520) by NSFC and Early Career Scheme (Project number: 27301921), RGC, Hong Kong.}}
\maketitle
\begin{abstract}
Computing the dispersion relation for two-dimensional photonic crystals is a notoriously
challenging task: It involves solving parameterized Helmholtz eigenvalue problems with
high-contrast coefficients. 
To resolve the challenge, we propose a novel $hp$-adaptive sampling scheme that can detect
singular points via adaptive mesh refinement in the parameter domain, and meanwhile, allow for
adaptively enriching the local polynomial spaces on the elements that do not contain singular
points. In this way, we obtain an element-wise interpolation on an adaptive mesh. We derive
an exponential convergence rate when the number of singular points is finite, and a
first-order convergence rate otherwise. Numerical tests are provided to illustrate its performance.
\end{abstract}

\section{Introduction}
This work is concerned with propagation of electromagnetic waves within photonic crystals (PhCs). PhCs are
dielectric materials with a period size comparable to the wavelength \cite{joannopoulos2008molding}. The
properties of waves in PhCs depend heavily on their frequencies, and for certain frequencies, they may be
unable to propagate through PhCs, leading to the band gap phenomenon. This intriguing feature
has led to the development of important applications, e.g., optical transistors, photonic fibers, and low-loss
optical mirrors \cite{yanik2003all,russell2003photonic,wang2023analytical,labilloy1997demonstration}.
We focus on two-dimensional (2D) PhCs, which are homogeneous along the $z$ axis and have high-contrast
dielectric columns or holes in dielectric materials within the $x$-$y$ plane.

The propagation of electromagnetic waves is governed by the Maxwell system. For 2D infinite
crystals with perfect periodicity, Bloch theorem states that the problem is reduced to 2D parameterized
Helmholtz eigenvalue problems defined over the unit cell with periodic boundary conditions. The
associated parameter is the so-called wave vector $\mathbf{k}$, which varies in the irreducible
Brillouin zone (IBZ) \cite{kuchment1993floquet}. The $n$th band function $\omega_n$ (or the
square root of the $n$th largest eigenvalue up to a constant) is a function of the wave vector for all $n\in\mathbb{N}_+$.
The band gap is the distance between two adjacent band functions. Thus, computing the $n$th band
function $\omega_n(\mathbf{k})$ involves solving infinite many Helmholtz eigenvalue problems defined on
the unit cell. To obtain the band gap, it is necessary for the permittivity to take different values
in the inclusion and background of the unit cell, and moreover, the ratio between these values, known as
{contrast}, should be large. However, solving Helmholtz eigenvalue problems with high-contrast and
piecewise constant coefficients poses significant challenges. Numerically, one possible approach is to
limit the parameter $\mathbf{k}$ to a coarse mesh of the IBZ or even its edges only. However, it is
inadequate for identifying small band gaps, e.g., avoided crossings, which are important in certain
application areas of photonic materials. For example, it is important to avoid the smallest band gaps
when manufacturing large-pitch photonic crystal fibers
\cite{eidam2011fiber}, whereas leakage channel fibers exploit the smallest band gaps \cite{dong2009extending}.
While directly solving the eigenvalue problem over a very fine parametric discretization is failure-safe,
it is highly inefficient. Band structure diagram paths based on crystallography categorize
wave vectors by their symmetry and provide preferred band paths for band function reconstruction
\cite{hinuma2017band}. These important wave vectors are located at the face and edge centers, and
vertices of Brillouin zone. However, there is no guaranteed accuracy. 

The main goal of this paper is to provide a numerical method with guaranteed fast convergence rate. 
More specifically, let $\mathcal{B}$ be the parameter domain, i.e., IBZ. Let $\omega(\mathbf{k})\in C(\mathcal{B})$ be a band function
and $P_n(\mathcal{B})$ the space of global polynomial functions of degree at most $n$. The
best approximation to $\omega(\mathbf{k})$ in $P_n(\mathcal{B})$ with respect to the uniform norm
$\Vert\cdot\Vert_{\infty,\mathcal{B}}$ on $\mathcal{B}$ can be measured by
\begin{equation*}
    E_n(\omega):=\inf_{p\in P_n(\mathcal{B})}\Vert \omega-p\Vert_{\infty,\mathcal{B}}.
\end{equation*}
If $\omega\notin C^1(\mathcal{B})$, then $E_n(\omega)$ cannot decay
faster than first order with respect to $n$, i.e., $E_n(\omega)=\mathcal{O}(n^{-1})$
\cite[Chapter 7]{timan2014theory}. Furthermore, this decay rate may not be improved in a subdomain
where the band function is analytic \cite{kadec1963distribution}. The decay rate is
only $\mathcal{O}(N^{-1/2})$ in terms of the number $N$ of sampling points $N$. In
\cite{wang2023dispersion}, we have developed a global polynomial interpolation method (GPI)
to reconstruct the band function, focusing on selecting sampling points to achieve this slow convergence rate.
Thus, there is an imperative need to develop more effective computational techniques for
band structure reconstruction.

To improve the convergence of GPI, it is critical to exploit piecewise analyticity
of band functions and to utilize local polynomial interpolation. This work focuses
on adaptively selecting the wave vectors. Inspired by the $h$,
$p$, and $hp$ versions of FEM \cite{babuvska1994p}, one promising approach is to discretize the
parameter domain $\mathcal{B}$ with a triangular mesh and approximate the band functions
locally on each element. In the $h$ version of local interpolation, the mesh sizes tend
to zero uniformly and the degree of the polynomials is fixed at a small
number, and the convergence rate is still $\mathcal{O}(N^{-1/2})$. In the $p$ version
of local interpolation, the mesh is fixed and the degree of polynomials tends to infinity,
and it can also achieve a convergence rate $\mathcal{O}(N^{-1/2})$. Our approach combines
these two strategies and consists of
two steps. First, we generate an adaptive mesh of $\mathcal{B}$ by Algorithm \ref{alg1}
such that elements containing singular points are refined at every iteration. Second, given a
slope parameter $\mu$, we assign an element-wise polynomial degree according to the number
of refinement of each element, and develop conforming element-wise interpolation in
Algorithm \ref{alg2}. To analyze the convergence of the method, we first provide in
Theorem \ref{property} the condition on the slope parameter $\mu$ such that the element-wise
interpolation error in the element with singular points dominates. Then we derive an
exponential convergence rate in Theorem \ref{exp} when the number of singular points is finite
and a first-order convergence rate in Theorem \ref{algebraic} otherwise.

Overall, the proposed $hp$-adaptive sampling algorithm significantly improves the computational
efficiency and accuracy of band function approximation in PhCs. The adaptive refinement of the
parameter domain and element-wise interpolation can effectively capture singularities and
improve approximation accuracy, while the mechanism for selecting polynomial degrees can
further enhance the convergence rate. Additionally, the computation of band functions at each selected
wave vector is independent and embarrassingly parallel. This work provides
a novel approach for photonic crystal band function computation, which holds
potential for designing and optimizing photonic crystal structures.

The rest of the paper is organized as follows. In Section \ref{sec:model}, we recap
the derivation of band functions for 2D PhCs and regularity.
We propose in Section \ref{sec:hp-sampling} an $hp$-adaptive sampling algorithm, which
adaptively refines elements containing singular points and allocates proper polynomial
degrees for local Lagrange polynomial interpolation. We define the element-wise interpolation
over the adaptive mesh and with a given element-wise degree in Section \ref{sec:elementwise-interp}.
In Section \ref{sec:convergence}, we present a convergence analysis of the algorithm.
We present extensive numerical experiments in Section \ref{sec:simulation} to
complement the analysis. Finally, we give concluding remarks in Section \ref{sec:conclusion}.
Throughout, the notation $A \lesssim B$ stands for that $A/B$ is bounded above by a constant
independent of the number of sampling points $N$ and the mesh size $h$.

\section{Problem formulation}\label{sec:model}
\subsection{2D Band function computation}
First we briefly derive  the parameterized Helmholtz eigenvalue problem from the Maxwell equations, which was formulated in \cite{wang2023dispersion} to estimate band functions. For further details, we refer readers to  \cite{bao2001mathematical,wang2023dispersion}. In the SI convention, the time harmonic Maxwell equations for linear, non-dispersive, and nonmagnetic media, with free charges and free currents, consist of a system of four equations  \cite{jackson1999classical},
\begin{subequations}
\begin{align}
    \nabla\times\mathbf{E}(\mathbf{x})-i\omega\mu_{0}\mathbf{H}(\mathbf{x})=0,\label{harm1}    \\
    \nabla\times\mathbf{H}(\mathbf{x})+i\omega\epsilon_{0}\epsilon(\mathbf{x})\mathbf{E}(\mathbf{x})=0 , \label{harm2}\\  \nabla\cdot\left(\epsilon(\mathbf{x})\mathbf{E}(\mathbf{x})\right)=0,\label{harm3}\\
    \nabla\cdot\mathbf{H}(\mathbf{x})=0.\label{harm4}
\end{align}
\end{subequations}
Here $\mathbf{x}\in\mathbb{R}^3$, $\mathbf{E}$ is the electric field, $\mathbf{H}$ the magnetic field and $\mathbf{D}$ the electric displacement field. The scalar $\omega\geq 0$ is the frequency of the electromagnetic wave, $\mu_0$ the vacuum permeability, $\epsilon_0$ the vacuum permittivity, and $\epsilon\in L^{\infty}(\mathbb{R}^3;\mathbb{R}^+)$ the relative permittivity.

Applying the curl operator to \eqref{harm1} and using \eqref{harm2}, we obtain
\begin{equation}\label{E}
 \nabla\times\left(\nabla\times\mathbf{E}(\mathbf{x})\right)-(\omega c^{-1})^2\epsilon(\mathbf{x})\mathbf{E}(\mathbf{x})=0,
\end{equation}
with $\epsilon_{0}\mu_{0}=c^{-2}$, where $c$ is the speed of light. Similarly, applying 
the curl operator to \eqref{harm2} and using \eqref{harm1}, we obtain
\begin{equation}\label{H}
    \nabla\times\left(\epsilon(\mathbf{x})^{-1}\nabla\times\mathbf{H}(\mathbf{x})\right)-\left(\omega c^{-1}\right)^2\mathbf{H}(\mathbf{x})=0.
\end{equation}

\begin{figure}[hbt!]
\centering
\subfigure[PhCs with square lattice]{\label{square material}
\includegraphics[width=.35\textwidth]{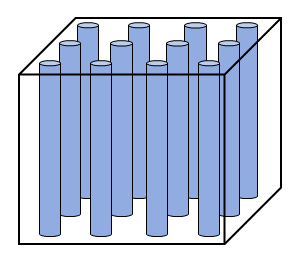}
}%
\subfigure[PhCs with hexagonal lattice]{\label{hex material}
\includegraphics[width=.45\textwidth]{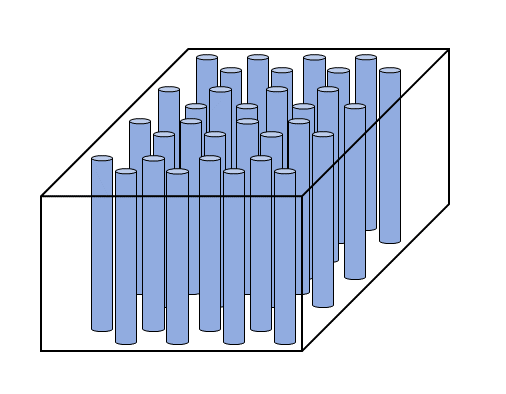}
}%
\centering
\caption{A schematic illustration of 2D PhCs. The materials have dielectric columns homogeneous along the
$z$ direction and periodic along $x$ and $y$ directions.}\label{material}
\end{figure}

In this work, we focus on 2D PhCs. In  Figure \ref{material}, we  show 2D PhCs with a periodic arrangement of dielectric columns within a dielectric material. In practice, both have finite extensions in the $z$ direction and finite periodicities in the $x$-$y$ plane. Nevertheless, it is commonly assumed that the material extends infinitely in the plane perpendicular to the columns. Thus, in the setting of 2D PhCs in Figure \ref{material}, the relative permittivity $\epsilon(\mathbf{x})$ is assumed to be independent of the $z$ direction. Then we can split the electromagnetic fields $\mathbf{E}=(E_1,E_2,E_3)$ and $\mathbf{H}=(H_1,H_2,H_3)$ in \eqref{E} and \eqref{H} into the transverse electric (TE) mode with $H_1=H_2=E_3=0$ and transverse magnetic (TM) mode with $E_1=E_2=H_3=0$. Each mode is a scalar eigenvalue problem,
\begin{align}
    -\nabla \cdot (\epsilon(\mathbf{x})^{-1}\nabla H(\mathbf{x}))-(\omega c^{-1})^2 H(\mathbf{x})&=0,\quad \mathbf{x}
\in\mathbb{R}^2,&&\textbf{(TE mode)},\label{TE}\\
    -\Delta E(\mathbf{x})-(\omega c^{-1})^2 \epsilon(\mathbf{x})E(\mathbf{x})&=0,\quad \mathbf{x}
\in\mathbb{R}^2,&& \textbf{(TM mode)}.\label{TM}
\end{align}

\begin{figure}[hbt!]
\centering
\subfigure{\label{lattice(1)}
\includegraphics[width = .3\textwidth,trim={14cm 4.5cm 12cm 4.5cm},clip]{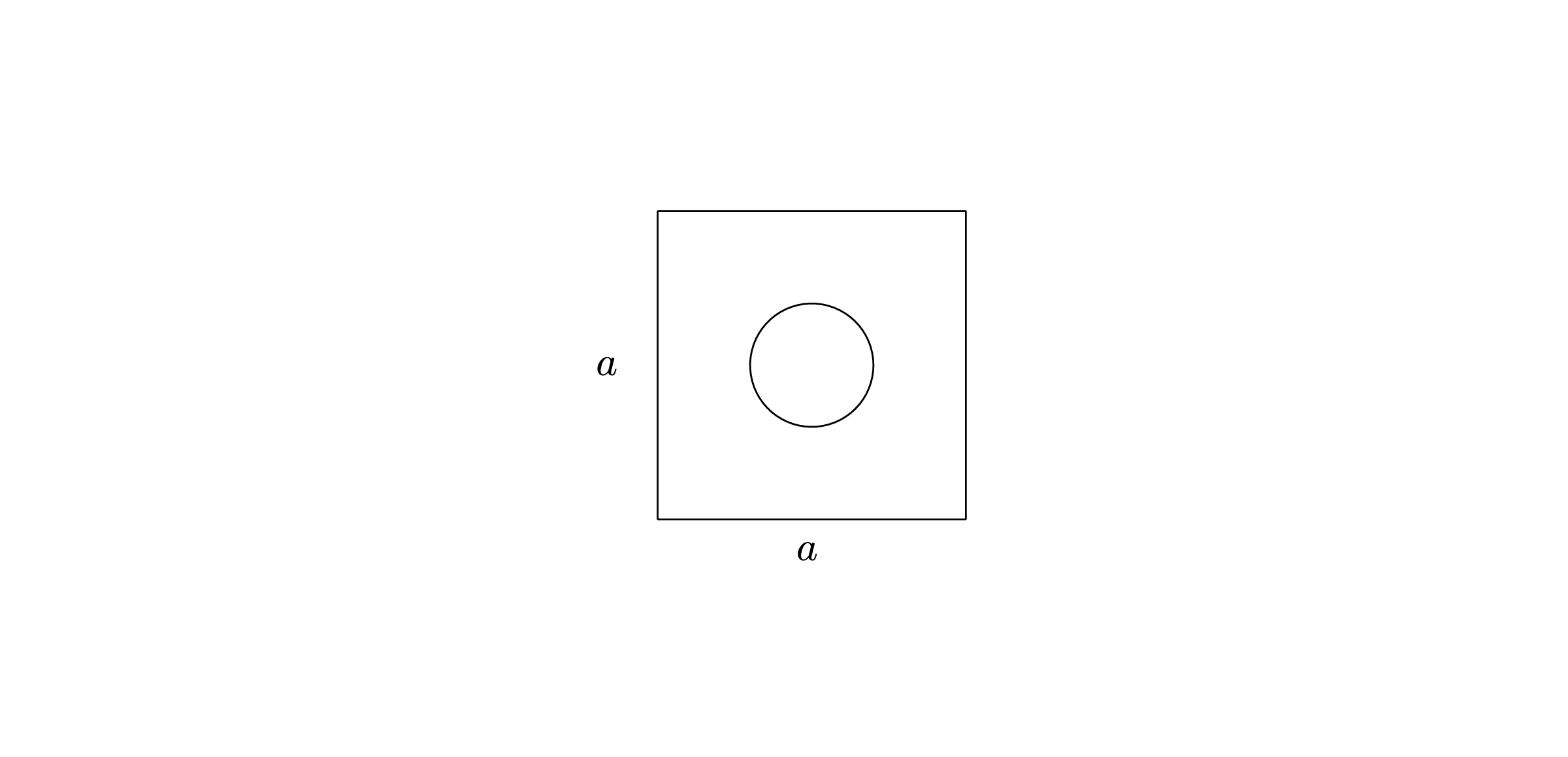}
}%
\subfigure{\label{lattice(2)}
\includegraphics[width = .3\textwidth,trim={14cm 4.5cm 12cm 4.5cm},clip]{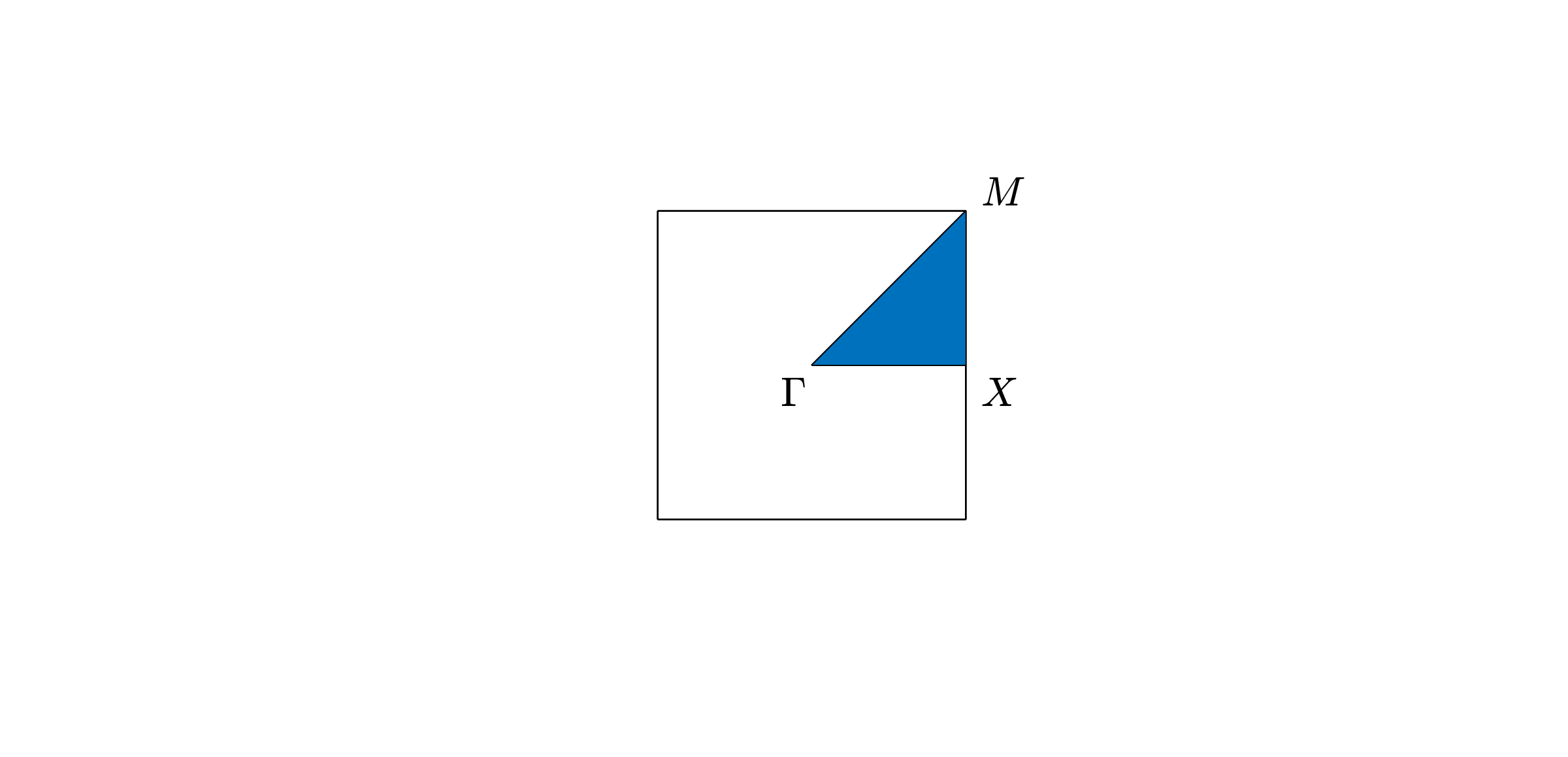}
}%
\centering
\caption{Square lattice: $\Omega$ (left) and the corresponding Brillouin zone (right). The IBZ $\mathcal{B}$} is the shaded triangle with vertices $\Gamma=(0,0)$, $X=\frac{1}{a}(\pi,0)$ and $M=\frac{1}{a}(\pi,\pi)$.\label{lattice}
\end{figure}

\begin{figure}[hbt!]
\centering
\subfigure{\label{lattice2(1)}
\includegraphics[width = .27\textwidth,trim={14cm 4.5cm 13cm 4.5cm},clip]{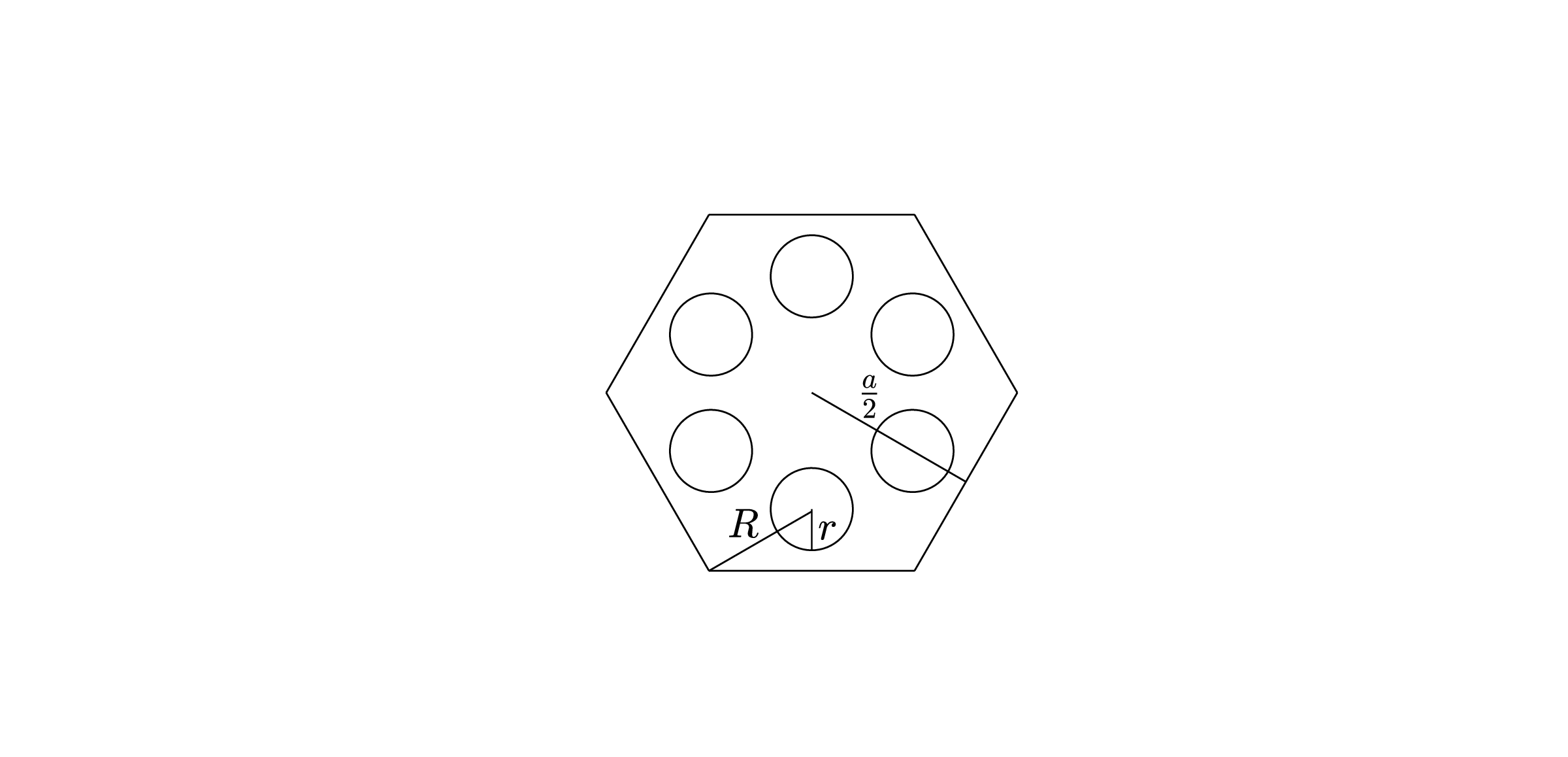}
}%
\subfigure{\label{lattice2(2)}
\includegraphics[width = .27\textwidth,trim={14cm 4.5cm 12cm 4.5cm},clip]{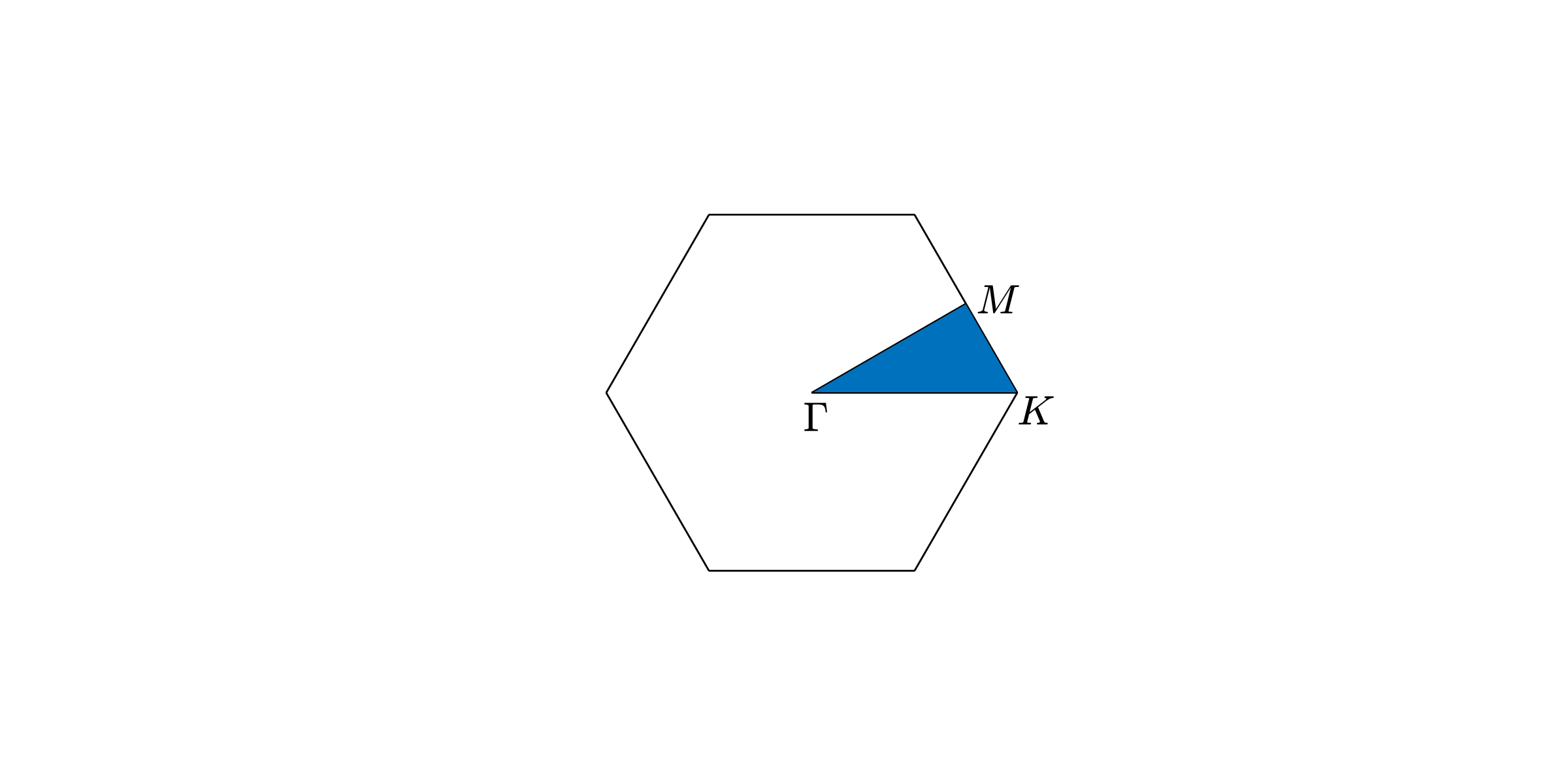}
}%
\centering
\caption{Hexagonal lattice: $\Omega$ (left) and the corresponding Brillouin zone (right). The IBZ $\mathcal{B}$} is the shaded triangle with vertices $\Gamma=(0,0)$, $K=\frac{1}{a}(\frac{4}{3}\pi,0)$ and $M=\frac{1}{a}(\pi,\frac{\sqrt{3}}{3}\pi)$.\label{lattice3}
\end{figure}

2D PhCs possess a discrete translational symmetry in the $x$-$y$ plane
(cf.  Figure \ref{material}) \cite{joannopoulos2008molding}, i.e., the relative
permittivity $\epsilon(\mathbf{x})$ satisfies
\begin{equation*}
   \epsilon(\mathbf{x}+c_1\mathbf{a}_1+c_2\mathbf{a}_2)=\epsilon(\mathbf{x}),\quad
    \forall\mathbf{x}\in\mathbb{R}^2\text{ and }c_1,c_2\in\mathbb{Z}.
\end{equation*}
The primitive lattice vectors, denoted by $\mathbf{a}_1$ and $\mathbf{a}_2$, are the shortest possible vectors that fulfill this condition and they span the fundamental periodicity domain $\Omega$, also known as unit cell, cf. Figures \ref{lattice} and \ref{lattice3}. In Figure \ref{lattice}, primitive lattice vectors of the square lattice are $\mathbf{a}_i = a{e}_i$ for $i=1,2$, and in Figure \ref{lattice3}, that of the hexagonal lattice $\mathbf{a}_1 = \frac{a}{4}({e}_1+\sqrt{3}{e}_2)$ and $\mathbf{a}_2 = \frac{a}{4}({e}_1-\sqrt{3}{e}_2)$. Here, $(e_i)_{i=1,2}$ is the canonical basis in $\mathbb{R}^2$ and $a\in \mathbb{R}^+$ lattice constant. The reciprocal lattice vectors are defined by
\begin{equation}
  \mathbf{b}_i \cdot \mathbf{a}_j=2\pi\delta_{ij},\quad \text{ for }i,j=1,2,
\end{equation}
which generate the so-called reciprocal lattice. The elementary cell of the reciprocal lattice is the (first) Brillouin zone $\mathcal{B}_F$, i.e., the region closer to a certain lattice point than to any other lattice points in the reciprocal lattice.

Bloch's theorem \cite{kittel2018introduction} states that in periodic crystals, wave functions take the form of a plane wave modulated by a periodic function. Thus they can be written as
$\Psi(\mathbf{x})=e^{i\mathbf{k}\cdot\mathbf{x}}u(\mathbf{x})$,
where $\Psi$ is the wave function, $u(\mathbf{x})$ is periodic sharing the periodicity of the crystal lattice and $\mathbf{k}$ is the wave vector varying in the  Brillouin zone $\mathcal{B}_F$. The periodic condition of $u(\mathbf{x})$ implies that each wave function can be determined by its values within the unit cell $\Omega$. Thus, the solutions to \eqref{TE} and \eqref{TM} can be expressed as
$H(\mathbf{x})=e^{i\mathbf{k}\cdot \mathbf{x}}u_1(\mathbf{x})$ and $E(\mathbf{x}) =e^{i\mathbf{k}\cdot \mathbf{x}}u_2(\mathbf{x})$ for some periodic functions $u_1(\mathbf{x})$ and $u_2(\mathbf{x})$ in the unit cell $\Omega$, and the parameterized Helmholtz eigenvalue problems \eqref{TE} and \eqref{TM} reduce to
\begin{subequations}
\begin{align}
-(\nabla+i\mathbf{k})\cdot\left(\epsilon(\mathbf{x})^{-1}(\nabla+i\mathbf{k}) u_1(\mathbf{x})\right)-\left( \omega c^{-1}\right)^{2} u_1(\mathbf{x})&=0,\quad \mathbf{x}\in\Omega && \textbf{(TE mode)},\label{TE2}\\
   -(\nabla+i\mathbf{k})\cdot\left((\nabla+i\mathbf{k}) u_2(\mathbf{x})\right)-\left(\omega c^{-1}\right)^{2}\epsilon(\mathbf{x})u_2(\mathbf{x})&=0,\quad \mathbf{x}\in\Omega&& \textbf{(TM mode)},\label{TM2}
\end{align}
\end{subequations}
where $\mathbf{x}\in\Omega\subset\mathbb{R}^2$, $\mathbf{k}$ varies in the Brillouin zone (BZ), and $u_i(\mathbf{x})$ satisfies
the periodic boundary conditions $u_i(\mathbf{x})=u_i(\mathbf{x}+\mathbf{a_j})$ with $\mathbf{a_j}$ being the primitive lattice vector for $i,j=1,2$.
If the materials in the unit cell have additional symmetry, e.g. mirror symmetry, we can further restrict $\mathbf{k}$ to the triangular irreducible Brillouin zone (IBZ), denoted by $\mathcal{B}$. The BZs and IBZs of the square lattice and hexagonal lattice are shown in Figures \ref{lattice} and \ref{lattice3}.

In sum, we can formulate both parameterized Helmholtz problems \eqref{TE2} and \eqref{TM2} by
\begin{equation}\label{both}
    -(\nabla+i\mathbf{k})\cdot \alpha(\mathbf{x})(\nabla+i\mathbf{k})u(\mathbf{x})-\lambda\beta(\mathbf{x})u(\mathbf{x})=0, \quad \mathbf{x}\in\Omega,
\end{equation}
with $\Omega\subset\mathbb{R}^2$, $\mathbf{k} \in \mathcal{B}$, and $\lambda =(\omega c^{-1})^2$. In the TE mode, $u$ describes the magnetic field $H$ in $z$-direction and the coefficients $\alpha(\mathbf{x})$ and $\beta(\mathbf{x})$ are $ \alpha(\mathbf{x}):=\epsilon(\mathbf{x})^{-1}$ and $\beta(\mathbf{x}):=1$.
Similarly, in the TM mode, $u$ describes the electric field $E$ in the $z$-direction and the coefficients $\alpha(\mathbf{x})$ and $\beta(\mathbf{x})$ are
$ \alpha(\mathbf{x}):=1$ and $\beta(\mathbf{x}):=\epsilon(\mathbf{x})$.
The variational formulation of \eqref{both} reads: find a non-trivial eigenpair $(\lambda,u) \in (\mathbb{R},H^1_{\pi}(\Omega))$ for $\mathbf{k}\in\mathcal{B}$ such that
\begin{equation}\label{variational}
\left\{
\begin{aligned}   \int_{\Omega}\alpha(\nabla+i\mathbf{k})u\cdot(\nabla-i\mathbf{k})\Bar{v}-\lambda\beta u\Bar{v}\mathrm{d}x&=0, \quad \forall v \in H^1_{\pi}(\Omega)\\
\|u\|_{L^2_{\beta}(\Omega)}&=1.
    \end{aligned}
\right.
\end{equation}
Here, we define $L^2_{\beta}(\Omega)$ as the space of weighted square integrable
functions equipped with the norm
$\|f\|_{L^2_{\beta}(\Omega)}:=(\int_{\Omega} \beta(\mathbf{x})|f(\mathbf{x})|^2\mathrm{d}\mathbf{x})^{\frac{1}{2}}$.
Let $H^1(\Omega)\subset L^2_{\beta}(\Omega)$ be with square integrable gradient, equipped with the standard $H^1(\Omega)$-norm and $H^1_{\pi}(\Omega)\subset H^1(\Omega)$ is composed of functions with periodic boundary conditions.
Using the sesquilinear forms
\begin{align*}
    a(u,v)&:=\int_{\Omega}\alpha (\nabla+i\mathbf{k})u\cdot(\nabla-i\mathbf{k})\Bar{v}\,\mathrm{d}x\quad \mbox{and}\quad
    b(u,v):=\int_{\Omega}\beta u\Bar{v}\,\mathrm{d}x,
\end{align*}
problem \eqref{variational} is equivalent to finding a non-trivial eigenpair $(\lambda,u) \in (\mathbb{R},H^1_{\pi}(\Omega))$ for $\mathbf{k}\in\mathcal{B}$ such that
\begin{equation}\label{simply}
    \left\{
\begin{aligned}
a(u,v)&=\lambda b(u,v), \text{ for all }v \in H^1_{\pi}(\Omega) \\
\|u\|_{L^2_{\beta}(\Omega)}&=1.
\end{aligned}
\right.
\end{equation}
\subsection{Regularity of band functions in the parameter domain}
The regularity of band functions plays a crucial role in the convergence analysis. We
have the following result taken from \cite{wang2023dispersion}. Here, $ Lip(\mathcal{B})$
denotes the space of Lipschitz continuous functions in the IBZ $\mathcal{B}$ and
$\mathring{A}(\mathcal{B})$ the space of piecewise analytic functions with singular point set
of zero Lebesgue measure.
\begin{theorem}
\label{thm:lipchitz}
For 2D periodic PhCs, $\omega_n(\mathbf{k})\in Lip(\mathcal{B})\cap \mathring{A}(\mathcal{B})$ for all $n\in\mathbb{N}_+$: each band function $\omega_n(\mathbf{k})$ is analytic in $\mathcal{B}\backslash\mathbf{S}$, where $\mathbf{S}\subset \mathcal{B}$ consists of the origin and points of degeneracy (i.e., the set of points at which certain band functions intersect, with zero Lebesgue measure).
\end{theorem}

\begin{remark}\label{singular at 0}
Note that the singularity at $\mathbf{k}=\mathbf{0}$ differs from that of others. The latter are known as branch points while the former arises from the composition $g\circ \lambda$ with $g(x):=\sqrt{x}$. The derivative of the first band function $\omega_1(\mathbf{k})=c\cdot \sqrt{\lambda_1(\mathbf{k})}$, which takes the value $0$ at the origin, fail to be defined at $\mathbf{k}=\mathbf{0}$. To remove the effect of this singularity, we approximate $\lambda_j(\mathbf{k})$ for $j=1,\cdots,L-1$ with $1<L\in\mathbb{N}$, instead.
\end{remark}

\section{$hp$-adaptive sampling algorithm in the parameter domain}\label{sec:hp-sampling}
Due to the presence of singularities, the convergence of numerical methods based on global
or fixed-order local polynomial interpolation on uniform meshes may be slow.
To optimally deal with singular points of  band functions in PhCs, we propose a
novel efficient numerical algorithm using an $hp$-adaptive interpolation method. It is inspired by the celebrated $hp$
FEM \cite{babuvska1994p,ainsworth1997posteriori}, and employs nonuniformly distributed
polynomial degrees over graded, nonuniform meshes to discretize the problem.
We shall design an $hp$ interpolation strategy based on an element-wise quantity in order to
approximate the first $L-1$ band functions for $1<L\in\mathbb{N}$.

\subsection{Adaptive mesh generation}\label{Mesh design}
We first recall the concept of a regular mesh without hanging nodes \cite{schwab1998p}.
Let $D\subset\mathbb{R}^2$ be a polygon with straight sides. A mesh $\mathcal{T}=\{T_i\}_{i=1}^M$
is a partition of $D$ into open disjoint triangles $T_i$ with $\overline{D}=
\cup_{i=1}^M\overline{T_i}$. For each $T\in \mathcal{T}$, we denote by
$h_T$ its diameter and $\rho_T$ the diameter of the largest inscribed ball.
\begin{definition}
\label{regular mesh}
A mesh $\mathcal{T}=\{T_i\}_{i=1}^M$ is said to be a regular mesh without hanging nodes
if the following conditions hold.
{\rm(i)} For $i\neq j$, $\overline{T_i}\cap\overline{T_j}$
is either empty or it consists of a vertex or an entire edge of $T_i$ and {\rm(ii)}
There exists a constant $C>0$ such that $h_T\leq C\rho_T$.
\end{definition}

We shall develop a strategy to generate a family of nested regular triangular meshes $\mathcal{T}_n$, $n=1,2,\cdots$,
based on an element-wise quantity \eqref{indicator1} and a criterion for refinement \eqref{tolerance}. We
need several notations for the mesh $\mathcal{T}_n$. The mesh size of $\mathcal{T}_n$ is denoted by
$h_n:=\max_{T\in\mathcal{T}_n}h_T$. The collection of all faces and vertices over $\mathcal{T}_n$ is 
denoted by $\mathcal{F}_n$ and $\mathcal{V}_n$, respectively. For any $T\in\mathcal{T}_n$,
$\mathcal{F}_T:=\partial T\cap \mathcal{F}_n$ denotes its faces, and $F_T^j$ for $j=1,2,3$ the $j$th face,
$\mathcal{V}_T:=\partial T\cap \mathcal{V}_n$ its vertices,
and $V_T^j$ for $j=1,2,3$ the $j$th vertex. For any $F\in\mathcal{F}_n$,
$\mathcal{T}_F:=\cup\{T\in\mathcal{T}_n:F\in\mathcal{F}_T\}$ denotes the collection
of elements having $F$ as an edge. Finally, $P_m(T)$ is the space of polynomial functions of degrees
at most $m\in\mathbb{N}_+$ on element $T$, and $P_m(F)$ the space of polynomial functions
of degrees at most $m\in\mathbb{N}_+$ on each face $F\in\mathcal{F}_n$.

Since each band function is a piecewise analytic function of the wave vector, with singular points at branch points and the origin, the mesh design strategy should identify elements in which adjacent bands are in close proximity, and then refine them suitably to minimize the impact of singularities on the approximation. To this end, we employ an element-wise indicator that is based on the distance between two adjacent band functions
\begin{equation}\label{indicator1}
\eta(T):=\min_{1\leq q \leq L-1}\left\|\omega_q-\omega_{q+1}\right\|_{\infty,\mathcal{V}_T} \text{ for }T\in \mathcal{T}_n.
\end{equation}
That is, it uses the adjacent values on the vertices of each triangle element to ascertain the regularity of the band
functions within each element. Furthermore, we introduce an element-wise tolerance
\begin{equation}\label{tolerance}
     \operatorname{tol}_1(T)=\kappa h_T \max_{1\leq q \leq L} \Vert\nabla\omega_q\Vert_{\infty,T\backslash\mathbf{S}},
 \end{equation}
where $\mathbf{S}$ is the set of singular points and $\kappa\geq 2\sqrt{2} $ is a constant. The band function $\omega_q(\mathbf{k})$ and the group velocity $\nabla\omega_q(\mathbf{k})$ used in $\eta(T)$ and $\operatorname{tol}_1(T)$ are computed using the conforming Galerkin FEM, cf. \eqref{simply} and \eqref{derivative}.

Give a mesh $\mathcal{T}_n$, the collection of marked elements is
\begin{align}\label{eq:mark}
\mathcal{M}_n:=\cup\left\{
T\in\mathcal{T}_n: \eta(T)\leq \operatorname{tol}_1(T)
\right\}.
\end{align}
We employ the iterated newest vertex bisection (NVB) \cite{funken2011efficient} to refine the elements in $\mathcal{M}_n$. It involves repeatedly bisecting the longest edge of the marked element and inserting a new vertex at the midpoint until the desired level of refinement is achieved. Then the resulting mesh is sufficiently refined in areas where the singularities may impact the approximation accuracy.

Iterated NVB preserves the mesh regularity during refinement \cite{traxler1997algorithm}. However, to achieve a set of regular meshes without any hanging nodes, it is often necessary to refine certain elements beyond the initial refinement. Thus, the number of elements to be refined, denoted as $\# \mathcal{T}_{n+1}-\# \mathcal{T}_{n}$, is typically larger than the number of marked elements $\# \mathcal{M}_n$. Nonetheless, the cumulative number of elements added in this manner does not inflate the total number of marked elements \cite[Theorem 2.4]{binev2004adaptive}.
\begin{theorem}\label{the number of elements}
Let $\mathcal{T}_1$ be a regular mesh in the sense of Definition \ref{regular mesh}. For $n\geq 1$ let $\{\mathcal{T}_n\}_{n\geq 1}$ be any sequence of the refinement of $\mathcal{T}_1$, where $\mathcal{T}_{n+1}$ is a regular mesh generated from $\mathcal{T}_{n}$ by the iterated NVB with a subset of the marked elements $\mathcal{M}_n\subset \mathcal{T}_n$. Then there exists a constant $C$ solely depending on $\mathcal{T}_1$ such that
    \begin{equation}
        \# \mathcal{T}_{n}-\# \mathcal{T}_{1}\leq C\sum_{i=1}^{n-1}\#\mathcal{M}_i,\quad \forall n\geq 1.
    \end{equation}
\end{theorem}

We give the proposed mesh refinement procedure in Algorithm \ref{alg1}. The parameters $\operatorname{tol}_2$
and $n_{\max}$ denote the smallest allowed element size and the maximum number of iterations, respectively.
\begin{algorithm}[H]
\caption{Adaptive mesh refinement}\label{alg1}
\begin{algorithmic}
\Require initial mesh $\mathcal{T}_1$, 
   tolerance $\operatorname{tol}_2$, and maximum loop number $n_{\max}$
\Ensure $\mathcal{T}_n$
\State $n \gets 1$
\While{$n \leq n_{\max}$}

\State Set $\mathcal{M}_n=\emptyset$
\For{$T\in\mathcal{T}_n$}
\If{$\eta(T)\leq \operatorname{tol}_1(T)$ and $h_T\geq \operatorname{tol}_2$}

\State $\mathcal{M}_n=\mathcal{M}_n\cup\{T\}$
\EndIf
\EndFor

Refine the collection of marked elements $\mathcal{M}_n$ by iterated NVB and adjust the mesh until there are no hanging nodes.

$n \gets n+1$
\EndWhile
\end{algorithmic}
\end{algorithm}

\subsection{Reliability of Algorithm \ref{alg1}}\label{Reliability}
For the mesh $\mathcal{T}_n$ at the $n$th iteration, let
\begin{align*}
\mathcal{T}_n^R:=\cup\left\{T\in \mathcal{T}_n:T\cap \mathbf{S}=\emptyset\right\}\text{ and }
\mathcal{T}_n^S:=\mathcal{T}_n\backslash \mathcal{T}_n^R
\end{align*}
be a partition of the mesh $\mathcal{T}_n$ based on the singular point set $\mathbf{S}$. The next result guarantees that Algorithm \ref{alg1} refines all the elements in $\mathcal{T}_n^S$ at the next iteration.
\begin{theorem}[$\mathcal{T}_n^S\subset \mathcal{M}_n$]
\label{refining mechanism}
 Algorithm \ref{alg1} with the element-wise quantity \eqref{indicator1} and the tolerance \eqref{tolerance} guarantees that all elements containing branch points are marked for refinement.
\end{theorem}
\begin{proof}
The proof relies on piecewise analyticity and Lipschitz continuity of band functions. Let $\mathcal{T}_n$ be a regular mesh generated by Algorithm \ref{alg1} for $n=1,\cdots$. For any band function $\omega(\mathbf{k})$ and any element $T\in \mathcal{T}_n^R$, by the mean value theorem, we obtain
\begin{equation}\label{inequality good}
    \max_{\substack{\mathbf{k}\in T, \widehat{\mathbf{k}}\in \mathcal{V}_T}}|\omega(\mathbf{k})-\omega(\widehat{\mathbf{k}})|\leq \sqrt{2}h_T\cdot \Vert\nabla\omega\Vert_{\infty,T}.
\end{equation}
If $T\in \mathcal{T}_n^S$, then the first derivative of band functions may be discontinuous on singular points.
The singular points on the line segment connecting $\mathbf{k}$ and $\widehat{\mathbf{k}}$ divide the segment
into several small intervals. By applying the mean value theorem on each subinterval, we obtain
\begin{equation}\label{inequality bad}
     \max_{\substack{\mathbf{k}\in T,\widehat{\mathbf{k}}\in \mathcal{V}_T}}|\omega(\mathbf{k})-\omega(\widehat{\mathbf{k}})|\leq \sqrt{2}h_T\cdot \Vert\nabla\omega\Vert_{\infty, T\backslash\mathbf{S}}.
\end{equation}
Note that $T\notin \mathcal{M}_n$ implies any two adjacent band functions $\omega_p(\mathbf{k})$ and $\omega_{p+1}(\mathbf{k})$ for $1\leq p<L-1$ satisfy
\begin{align*}
    \min_{\widehat{\mathbf{k}}\in \mathcal{V}_T}|\omega_p(\widehat{\mathbf{k}})-\omega_{p+1}(\widehat{\mathbf{k}})|> \operatorname{tol}_1(T)=\kappa h_T\cdot \max_{1\leq q \leq L} \Vert\nabla\omega_q\Vert_{\infty,T\backslash\mathbf{S}}.
\end{align*}
By the triangle inequality, this further implies for any $\widehat{\mathbf{k}}\in \mathcal{V}_T$ and $\mathbf{k}\in T$,
\begin{align*}
    \operatorname{tol}_1(T)&< |\omega_p(\widehat{\mathbf{k}})-\omega_{p+1}(\widehat{\mathbf{k}})|\\
    &\leq |\omega_p(\widehat{\mathbf{k}})-\omega_{p}(\mathbf{k})|+|\omega_p(\mathbf{k})-\omega_{p+1}(\mathbf{k})|+|\omega_{p+1}(\mathbf{k})-\omega_{p+1}(\widehat{\mathbf{k}})|.
\end{align*}
This, together with \eqref{inequality good} and \eqref{inequality bad}, leads to
\begin{align}
    \min_{\mathbf{k}\in T}|\omega_p(\mathbf{k})-\omega_{p+1}(\mathbf{k})|> (\kappa-2\sqrt{2}) h_T\cdot \max_{1\leq q \leq L} \Vert\nabla\omega_q\Vert_{\infty,T\backslash\mathbf{S}}\geq 0.
\end{align}
That is, there is no branch points in $T$ and $T\notin\mathcal{T}_n^S$.
Hence, we have proved that $T\notin\mathcal{T}_n^S$ if $T\notin\mathcal{M}_n$, or equivalently $\mathcal{T}_n^S\subset \mathcal{M}_n$.
\end{proof}

Computing the element-wise tolerance $\operatorname{tol}_1(T)$ \eqref{tolerance} is impractical since it involves estimating $\Vert\nabla\omega(\mathbf{k})\Vert_{\infty}$ for all $\mathbf{k}\in T\backslash\mathbf{S}$. However, numerical experiments
have shown that the slope of the first several band functions (a.k.a. group velocity) does not exhibit high
oscillations. Hermann et al \cite{hermann2001photonic} showed that in the first three bands, higher bands
have low values of the group velocity and all of them are curves without
significant fluctuations. Therefore, as a substitute for $\operatorname{tol}_1(T)$, we use a feasible tolerance
\begin{align}\label{eq:tol1-num}
\widehat{\operatorname{tol}}_1(T):=\kappa h_T\max_{1\leq q \leq L}\Vert\nabla\omega_q\Vert_{\infty,\mathcal{V}_T\backslash\mathbf{S}}.
\end{align}
Computing $\widehat{\operatorname{tol}}_1(T)$ only requires the band function and its partial derivative over the vertices. Moreover, the group velocity can be derived directly from the following expression \cite{wang2023dispersion}
\begin{equation}\label{derivative}
  \frac{\partial \omega}{\partial k_i}=\frac{2k_i m_\alpha(u,u)+m_{\alpha i}(u,u)}{2\omega },
\end{equation}
where the bilinear forms $m_{\alpha i}(\cdot,\cdot)$ and $m_{\alpha}(\cdot,\cdot)$ are defined as
\begin{align*}
&m_{\alpha i}(u,v)=\int_{\Omega}i\alpha\left(u
\frac{\partial\Bar{v}}{\partial x_i}-\Bar{v} \frac{\partial u}{\partial x_i}\right)\,\mathrm{d}\mathbf{x}, \quad i=1,2\quad \mbox{and}\quad
m_{\alpha}(u,v)=\int_{\Omega}\alpha u \Bar{v}\,\mathrm{d}\mathbf{x}.
\end{align*}

\subsection{Conforming finite element space}\label{Degree design}

Let $\{\mathcal{T}_n\}_{n=1}^{\infty}$ be a sequence of
nested regular triangular meshes generated by Algorithm \ref{alg1}. We shall construct an
element-wise polynomial space for the mesh $\mathcal{T}_n$. We first define the layer of
an element $T\in\mathcal{T}_n$ given by Algorithm \ref{alg1}.
\begin{definition}
\label{layer}
 The layer $\ell_T$ of an element $T\in \mathcal{T}_n$ is $\ell_T:=n-r$, where $0\leq r\leq n-1$ is the number of refinements that have been performed on this element in Algorithm \ref{alg1}.
\end{definition}

Note that the layer $\ell_T$ of an element $T\in\mathcal{T}_n$ implies the local regularity of the eigenfunction $\lambda_j(\mathbf{k})$ for $j=1,\cdots,L-1$. Since the singular points fall into the marked element $\mathcal{M}_n$, we employ polynomials of lower degrees for the local interpolation in $\mathcal{M}_n$. Note also that the first band function vanishes at the origin, and is non-negative, sharing the symmetry of the lattice. This leads to the expression of $\left(\omega_1(\mathbf{k})\right)^2$ similar to an elliptic cone. More precisely, since $\left(\omega_1(\mathbf{k})\right)^2$ is analytic at $\mathbf{k}=\mathbf{0}$ if $\omega_1(\mathbf{0})\neq \omega_2(\mathbf{0})$, its Taylor expansion at the origin is
\begin{equation*}
\left(\omega_1(\mathbf{k})\right)^2=\left(\omega_1(k_1,k_2)\right)^2=\sum_{i=1}^\infty\sum_{j=0}^ia_{ij}\left(k_1\right)^{2j}\left(k_2\right)^{2i-2j},
\end{equation*}
with some $a_{ij}\in\mathbb{R}$ for $i\in\mathbb{N}$ and $j=0,\cdots, i$. Hence $\lambda_1(\mathbf{k})=( \omega_1(\mathbf{k})c^{-1})^2$ is quadratic near the origin, and we use polynomials of minimal degree $2$ to interpolate eigenvalues over the marked elements $\mathcal{M}_n$.
These observations motivate the following definition of the element-wise space of polynomials.
The notation $\left\lceil \cdot \right\rceil$ denotes the ceiling function.
\begin{definition}\label{degree design strategy}
For any element $T\in\mathcal{T}_n$, we take $P_{n_T}(T)$ as the local interpolation space, with
$n_T:=\max\left\{2,\lceil \mu {\ell}_T \rceil  \right\}$  if $T\notin \mathcal{M}_n$, and
$n_T:=2$ otherwise, where the parameter $\mu>0$ is determined in \eqref{eq:mu-cond} in Section 
\ref{sec:convergence} and $\mathcal{M}_n$ is the collection of marked elements \eqref{eq:mark}.
\end{definition}

In this way, each element $T$ is assigned a total degree $n_{T}$. However, neighboring elements may have different degrees on the common face. Suitable constraints are required to ensure the inter-element continuity, for which it is essential to associate with each face in the mesh a face degree.
\begin{definition}[Facewise space of polynomial functions]\label{edge degree design strategy}
Let $P_{m_F}(F)$ be the local interpolation space over $F$ for $F\in\mathcal{F}_n$ with
$m_F:=\min_{T\in\mathcal{T}_F}\{n_T\}$.
This implies
$\max_{F\in\mathcal{F}_T}\{m_F\}\leq n_T, \text{
 for all }T\in \mathcal{T}_n$.
\end{definition}
Finally, on a given mesh $\mathcal{T}_n$ with element-wise and facewise degrees $\{n_T\}_{T\in\mathcal{T}_n}$ and $\{m_F\}_{F\in\mathcal{F}_n}$, we define the conforming FEM space
\begin{align}\label{def:vn}
V_n:=\left\{v\in C(\mathcal{B}): v|_{T}\in P_{n_T}(T) \text{ and } v|_{F}\in P_{m_F}(F),\, \forall T\in\mathcal{T}_n \text{ and }F\in \mathcal{F}_n\right\},
\end{align}
which is a space of continuous piecewise polynomial functions on $\mathcal{B}$. Let $N$ be its dimension. 
We aim at approximating the eigenvalues $\lambda$, for $j=1,\cdots,L-1$, over $V_n$ by element-wise Lagrange interpolation.

\section{Element-wise interpolation over the adaptive mesh}\label{element-wise interpolation}\label{sec:elementwise-interp}
We develop Lagrange interpolation over each element $T\in\mathcal{T}_n$. For each element, we first map its closure $\overline{T}$ into the reference equilateral triangle $\hat{T}=\{\mathbf{z}:= (x,y):0\leq y\leq (1+x)\sqrt{3},-1\leq x\leq 0\text{ or } 0\leq y\leq (1-x)\sqrt{3},0\leq x\leq 1\}$ and then derive Lagrange interpolation on $\hat T$.
The vertices of $\hat{T}$ are $\mathbf{\hat{z}}_1=(-1,0)$, $\mathbf{\hat{z}}_2=(1,0)$ and $\mathbf{\hat{z}}_3=(0,\sqrt{3})$, with the face opposite to $\mathbf{\hat{z}}_i$ denoted by $\hat{F}^i$. For any $\{m,m_1,m_2,m_3\}\in \mathbb{N}_{+}^4$ with $m\geq m_i$ for $i=1,2,3$, \eqref{def:vn} implies that
the local space of polynomial functions over $\hat{T}$ is
\begin{align}\label{eq:loc-space}
\hat{P}:=\{v\in P_m(\hat{T}): v|_{\hat{F}^i}\in P_{m_i}(\hat{F}^i)\}.
\end{align}
Note that $\hat{P}\subset P_m(\hat{T})$ and $\hat{P}\subsetneq P_m(\hat{T})$ if $m_i<m$ for some $i=1,2,3$. Note also that the dimension of $P_{m}(\hat{T})$ is $\dim(P_{m}(\hat{T}))=(m+1)(m+2)/2$, and the dimension of $\hat{P}$ is
$\dim(\hat{P})=\dim(P_{m}(\hat{T}))-\sum_{i=1}^3(m-m_i)
=3+\big((m_1-1)+(m_2-1)+(m_3-1)\big)+(m-2)(m-1)/2.$
\subsection{Basis functions in the local polynomial space $\hat{P}$}
\label{Polynomial space}
We construct a set of basis functions for the local polynomial space $\hat{P}$, following \cite{schwab1998p}. We divide the basis functions into external and internal shape functions according to whether they vanish over the boundary $\partial\hat{T}$, and then divide the external shape functions further into nodal and side ones. The definition of nodal, side and internal shape functions is standard in the community of $p$-version and $hp$-adaptive FEMs, and it enables defining conforming FEM spaces with nonuniform degrees of polynomials over the mesh. We focus on developing an element-wise interpolation, and these shape functions play a role when formulating a special local polynomial space $\hat{P}$, which facilitates the definition of the element-wise interpolation in Section \ref{sec:loc-interpolation}.

First, there are three nodal shape functions on $\hat{T}$:
\begin{equation}\label{eq:nodal}
        N^1_1=(1-x-y/\sqrt{3})/2,\quad
        N^1_2=(1+x-y/\sqrt{3})/2\quad \mbox{and}\quad
        N^1_3=y/\sqrt{3}.
\end{equation}
Note that for $i=1,2,3$, we have
$N^1_i(\mathbf{\hat z}_i)=1$ and $N^1_i(\mathbf{z})=0$ for $\mathbf{z}\in \hat{F}^i$.
Second, there are $m_i-1$ side shape functions over $\hat{F}^i$, $i=1,2,3$, defined by
\begin{align}\label{eq:edgeB}
N^{m_i,j}_i&=N^1_{i+1}N^1_{i+2}\Phi_{j-2}(N^1_{i+2}-N^1_{i+1}),
\end{align}
with $\Phi_{j-2}(x)$ being a polynomial of degree $j-2$ for $j=2,\cdots,m_i$, with
the convention $N_{j+3}^1=N_j^1$ and $\mathbf{\hat{z}}_{j+3}=\mathbf{\hat{z}}_{j}$ for $j=1,2,3$. Note that $\text{supp}(N^{m_i,j}_i)$
is the interior of $\hat{F}^i$ and $\hat{T}$.
Then the set of external shape functions is given by
\begin{equation*}
    \mathcal{E}_m^{(m_1,m_2,m_3)}(\hat{T}):=\left\{N^1_j,N^{m_j,i}_j:\;1\leq j\leq 3\text{ and } 2\leq i\leq m_j\right\}.
\end{equation*}
By construction, the number of external basis functions is $3+\big((m_1-1)+(m_2-1)+(m_3-1)\big)$.

Next, we construct a set of internal basis functions that vanish over $\partial\hat{T}$. Let $b_{\hat{T}}:=N^1_1 N^1_2 N^1_3$ be the basic bubble function on the reference element $\hat{T}$ that vanishes over the edge $\partial\hat{T}$. Then the internal polynomial space $\mathcal{I}_m(\hat{T})$ is defined by
\begin{equation}\label{eq:internal-polynomial-space}
    \mathcal{I}_m(\hat{T}):=
    \text{span}\{b_{\hat{T}}v|v\in P_{m-3}(\hat{T})\}
    \text{ for }m\geq 3.
\end{equation}
The dimension of $\mathcal{I}_m(\hat{T})$ is $\hat{m}:={(m-1)(m-2)}/{2}$. Let one set of basis functions be denoted as
\begin{equation}\label{eq:internal-basis-functions}
\{\mathring{N}_i:i=1,\cdots, \hat{m}\}.
\end{equation}

\subsection{Local interpolation based on Lagrange polynomials}\label{sec:loc-interpolation}
Now following \cite{melenk2002hp}, we construct a suitable interpolation operator $\Pi: C(\hat{T}) \to \hat{P}$, such that
\begin{equation*}
\left\{
\begin{aligned}
\Pi f\big|_{\hat{F}^i}&\in P_{m_i}(\hat{F}^i)\text{ for } i=1,2,3,\\
\Pi f&\in \hat{P},
\end{aligned}
\right.\quad \forall f\in C(\hat{T}).
\end{equation*}
It consists of the external interpolation $E$ and internal interpolation $I$.

The external interpolation $E$ is an interpolation operator over Fekete points
along each face. Since the boundary points of Fekete points are the one-dimensional
Gauss-Lobatto points \cite{fejer1932bestimmung}, we get the Gauss-Lobatto interpolation
along each face $\hat{F}^i$ for $i=1,2,3$. Next, we extend the polynomial $Ef$ to the domain $\Omega$, and denote the extension by $\mathbf{E}(E f)$. With the shape functions in Section \ref{Polynomial space}, we can combine two procedures to obtain $\mathbf{E}(Ef)$. Specifically, we use the nodal and side shape functions to perform external interpolation along each face on Fekete points. For example, for $\hat{F}^i$ with $i=1,2, 3$, the nodal shape functions are $N_{i+1}^{1}$ and $N_{i+2}^{1}$, and we combine them with the $(m_i-1)$ side shape functions to obtain the interpolation on this face $\hat{F}^i$:
 \begin{equation}\label{eq:external-interp}
E f|_{\hat{F}^i}=\left(\sum_{j=2}^{m_i}c_{i,j}N_i^{m_i,j}+f(\mathbf{\hat{z}}_{i+1})N_{i+1}^1+f(\mathbf{\hat{z}}_{i+2})N_{i+2}^1\right)\Bigg|_{\hat{F}^i},
\end{equation}
where the coefficients $\{c_{i,j}\}_{j=2}^{m_i}$ are to be determined such that $E f|_{\hat{F}^i}$ interpolates
the $m_i$th order Fekete points on the face $\hat{F}^i$. 

The extension $\mathbf{E}(Ef)$ of $Ef$ is then given by
\begin{equation}\label{extension}
\mathbf{E}(E f)=\sum_{i=1}^{3}\sum_{j=2}^{m_i}c_{i,j}N_i^{m_i,j}+\sum_{i=1}^{3}f(\mathbf{\hat{z}}_{i})N_{i}^1.
\end{equation}
One can easily verify that $\mathbf{E}(E f)\big|_{\hat{F}^i}=Ef\big|_{\hat{F}^i}$ for $i=1,2,3$.

Third, we define the internal interpolation operator $I: C(\hat{T})\to \mathcal{I}_m(\hat{T})$ that interpolates the Fekete points $\{\mathbf{z}_j^{*}\}_{j=1}^{\hat{m}}$, using the set of internal basis functions $\{\mathring{N}_i:i=1,\cdots, \hat{m}\}$ \eqref{eq:internal-basis-functions}:
\begin{equation}\label{Vandermonde matrix}
 \begin{aligned}
\{\mathbf{z}_j^{*}\}_{j=1}^{\hat{m}}:=\arg\max_{\{\mathbf{z}_1,\mathbf{z}_2,\cdots,\mathbf{z}_{\hat{m}}\}\subset \hat{T}}\det(V(\mathbf{z}_1,\mathbf{z}_2,\cdots,\mathbf{z}_{\hat{m}})).
\end{aligned}
\end{equation}
Here, $V(\mathbf{z}_1,\mathbf{z}_2,\cdots,\mathbf{z}_{\hat{m}})\in \mathbb{R}^{\hat m\times \hat m}$ has entries
$V_{ij}=\mathring{N}_{j}(\mathbf{z}_i)$ for $i,j=1,\cdots,\hat{m}$.
Since $\mathring{N}_i$ vanishes on $\partial\hat{T}$ for $i=1,\cdots,\hat{m}$,
$\{\mathbf{z}_j^{*}\}_{j=1}^{\hat{m}}\subset\hat{T}\backslash\partial{\hat{T}}$. For any
$f\in C(\hat{T})$, the internal interpolation $I f$ is then given by
\begin{equation}\label{eq:inter-interp}
    I f=\sum_{i=1}^{\hat{m}}a_i\mathring{N}_i \text{ s.t. }
    If(\mathbf{z}_j^*)=f(\mathbf{z}_j^*) \text{ for all } j=1,\cdots, \hat{m}.
\end{equation}
Since the internal shape functions vanish along all three faces, $If|_{\partial\hat{T}}=0$ for all $f\in C(\hat{T})$.

Finally, we can define local interpolation operator for any function $f\in C(\hat{T})$ by,
\begin{equation}\label{eq:pi-n}
    \Pi f:=\mathbf{E}(E f)+I (f-\mathbf{E}(E f)).
\end{equation}
By construction, we can derive the following properties of $\Pi$.
\begin{proposition}
The element-wise interpolation $\Pi: C(\hat{T})\to \hat{P}$ satisfies the following properties.
{\rm(i)} $\Pi p=p$, for all $p\in \hat{P}$;
{\rm(ii)} $\Pi f\big|_{\hat{F}^i}=Ef\in P_{m_i}(\hat{F}^i)$ is the Lagrange interpolation on Fekete points along the face $\hat{F}^i$, for $i=1,2,3$.
\end{proposition}
\subsection{Element-wise interpolation}
For a given mesh $\mathcal{T}_n$, we define the operator $\Pi_{\mathcal{T}_n}: C(\mathcal{B})\to V_n$ by applying $\Pi$ to each element $T\in\mathcal{T}_n$:
\begin{equation}\label{eq:interp-global}
\Pi_{\mathcal{T}_n}f\big|_{T}:=\left(\Pi(f\circ M_T)\right)\circ M_T^{-1},\quad \forall T\in \mathcal{T}_n,
\end{equation}
where the affine transformation $M_T:=A_T\mathbf{z}+\mathbf{b}$ ($A_T\in\mathbb{R}^{2\times 2}$ is invertible and $\mathbf{b}\in \mathbb{R}^2$) maps the reference triangle $\hat{T}$ to the "physical" element $\Bar{T}$.

Now we can give the following algorithm for approximating band functions. By sequentially generating the adaptive mesh and the polynomial degree required for each element and each face, we can then use the interpolation formula \eqref{eq:interp-global}  to adaptively interpolate $\lambda_j$ and then generate the corresponding band function $\omega_j$ for $j=1,2,\cdots,L-1$.

\begin{algorithm}[H]
\caption{$hp$-adaptive sampling algorithm}\label{alg2}
\begin{algorithmic}
\Require $\mathcal{T}_n$: an adaptive mesh from Algorithm \ref{alg1};

    $\{\ell_T\}_{T\in \mathcal{T}_n}$: layer of each element defined in Definition \ref{layer};

    $\mu$: a positive parameter defined in Definition \ref{degree design strategy} satisfying condition \eqref{eq:mu-cond} to be introduced;
\Ensure $\omega_j$, $j=1,2,\cdots,L-1$

\State\hspace{\algorithmicindent}\hspace{0.2em} Choose the element-wise degree $\{n_T\}_{T\in \mathcal{T}_n}$ by Definition \ref{degree design strategy};

Choose the facewise degree $\{m_F\}_{F\in \mathcal{F}_n}$  by Definition \ref{edge degree design strategy};

Compute the local interpolation $\Pi\lambda_j$ by \eqref{eq:pi-n};

Compute the global interpolation $\Pi_{\mathcal{T}_n}\lambda_j$ by \eqref{eq:interp-global};

$\omega_j=c\cdot \sqrt{\lambda_j}$, $j=1,2,\cdots,L-1$.

\end{algorithmic}
\end{algorithm}

Now we compare the computational complexity of Algorithm \ref{alg2} with the global polynomial interpolation method \cite{wang2023dispersion} (ignoring their accuracy and computational cost of the adaptive mesh refinement in Algorithm \ref{alg1}). Let $\mathcal{T}_n$ and $V_n$ be the adaptive mesh generated by Algorithm \ref{alg1} and conforming FEM space \eqref{def:vn}, and let $N$ be the dimension of $V_n$, i.e., the number of sampling points. Suppose we use the same number $N$ of sampling points for both methods, and the computational complexity in solving the Helmholtz eigenvalue problem is the same. The main difference in the complexity arises from the Lagrange interpolation. To compute $N$ Lagrange polynomial functions using the classical Gauss elimination method involves $\mathcal{O}(N^{3})$ flops.
In contrast, Algorithm \ref{alg2} is based on element-wise interpolation, which has much lower computational complexity. Let the number of sampling points of each element $T\in\mathcal{T}_n$ be $N_T$. Then the local Lagrange interpolation involves $\mathcal{O}(N_T^{3})$ flops. Since there is no communication among each element, local Lagrange interpolation over each element can be performed in parallel. Hence, the total computational complexity is  $\max_{T\in\mathcal{T}_n}\mathcal{O}(N_T^{3})$ flops.

\section{Convergence analysis}\label{sec:convergence}
We establish in this section exponential and algebraic convergence of Algorithm \ref{alg2} in Theorems \ref{exp} and \ref{algebraic} when the number of iterations $n\to\infty$ in Algorithm \ref{alg1}, resting upon whether the number of crossings is finite or not. This indicates that Algorithm \ref{alg2} outperforms the global polynomial interpolation methods \cite{wang2023dispersion} with a convergence rate that is at least twice as fast. To establish Theorems \ref{exp} and \ref{algebraic}, we first derive the approximation property of the local polynomial interpolation operator $\Pi$ in Theorem \ref{local error theorem}. Then we formulate the element-wise approximation property of $\Pi_{\mathcal{T}_n}$ in Theorem \ref{property} and propose a proper assumption on the slope parameter $\mu$ \eqref{eq:mu-cond} such that the approximation error in the marked element patch $\mathcal{M}_n$ dominates.  

\subsection{Local approximation property}
Now we establish the approximation property of the local interpolation operator $\Pi$. Recall that for any $\{m,m_1,m_2,m_3\}\in \mathbb{N}_{+}^4$ with $m\geq m_i$ for $i=1,2,3$, the local space $\hat{P}$ is defined by \eqref{eq:loc-space}. First, we give the stability of the external interpolation and internal interpolation.
\begin{theorem}\label{theorem lebesgue}
Let $E$ and $I$ be the external interpolation and the internal interpolation defined by \eqref{eq:external-interp} and \eqref{eq:inter-interp}, respectively. Then the following stability estimates hold
\begin{align}
        \|E f\|_{\infty,\hat{F}^i}&\lesssim\log m_i\|f\|_{\infty,\hat{F}^i},&&\forall f\in C(\hat{F}^i),
        \label{eq:en}\\
        \|I f\|_{\infty,\hat{T}}&\leq \hat{m}\|f\|_{\infty,\hat{T}},&& \forall f\in C(\hat{T}). \label{eq:in}
\end{align}
\end{theorem}
\begin{proof}
Note that the Lebesgue constant of Gauss-Lobatto points in an interval (cf. \cite{sundermann1980lebesgue}) is $\mathcal{O}(\log m_i)$ and that of Fekete points in a triangle (cf. \cite{taylor2000algorithm}) has an upper bound $\hat{m}$. Then the desired result follows.
\end{proof}
Next, we derive the stability of the extension $\mathbf{E}$. To this end, we first introduce an inverse inequality for algebraic polynomials,  
\begin{lemma}
[Inverse estimate of Markov type \cite{wilhelmsen1974markov}]
\label{lemma ieq 1}
For all $m\in\mathbb{N}$ and $p\in P_m(\hat{T})$, there holds
\begin{equation*}
\Vert \nabla p\Vert_{\infty,\hat{T}}\leq 2m^2\Vert p\Vert_{\infty,\hat{T}}.
\end{equation*}
\end{lemma}
\begin{theorem}\label{theorem extension}
For each $p\in C(\partial \hat{T})$ with $p|_{\hat{F}^i}\in P_{m_i}(\hat{F}^i)$ for all $i=1,2,3$, the extension $\mathbf{E}$ defined in \eqref{extension} satisfies
    \begin{equation}\label{eq:ext}
        \|\mathbf{E}(p)\|_{\infty,\hat{T}}+m^{-2}\|\nabla \mathbf{E}(p)\|_{\infty,\hat{T}}\lesssim \|p\|_{\infty,\partial \hat{T}}.
    \end{equation}
Moreover, the extension map $\mathbf{E}:p\mapsto \mathbf{E}(p)$ is a bounded linear operator.
\end{theorem}
\begin{proof}
Since Lemma \ref{lemma ieq 1} implies $m^{-2}\|\nabla \mathbf{E}(p)\|_{\infty,\hat{T}}\leq 2\|\mathbf{E}(p)\|_{\infty,\hat{T}}$, it suffices to prove $\|\mathbf{E}(p)\|_{\infty,\hat{T}}\lesssim \|p\|_{\infty,\partial \hat{T}}$. We may assume $p$ vanishes at all vertices and all sides but one, e.g., the third face $\hat{F}^3=\{(x,0):-1\leq x\leq 1\}$. Then $\mathbf{E}(p)$ can be written in the form
$\mathbf{E}(p)=\sum_{i=2}^{m_3}c_{3,i}N_3^{m_3,i}$,
    where the unknown coefficients $\{c_{3,i}\}_{i=2}^{m_3}$ are to be determined such that
    \begin{align}\label{eq:inter1111}
    \mathbf{E}(p)=p \text{ on } \hat{F}^3.
    \end{align}
Since $N^1_{2}-N^1_{1}=x$, after plugging into \eqref{eq:edgeB}, we obtain
    \begin{equation*}
\mathbf{E}(p)=\sum_{i=2}^{m_3}c_{3,i}N^1_{1}N^1_{2}x^{i-2}=N^1_{1}N^1_{2}\sum_{i=2}^{m_3}c_{3,i}x^{i-2}.
\end{equation*}
Together with \eqref{eq:nodal} and \eqref{eq:inter1111}, this leads to
\begin{align*}
    \max_{(x,y)\in \hat{T}}|\mathbf{E}(p)|&=\max_{(x,y)\in \hat{T}}\left|N^1_{1}N^1_{2}\sum_{i=2}^{m_3}c_{3,i}x^{i-2}\right|\\
    &=\max_{\substack{0\leq x\leq 1,0\leq y\leq\sqrt{3}(1-x)\\-1\leq x\leq 0,0\leq y\leq\sqrt{3}(1+x)}}\frac{1}{4}\left(1-x-y/\sqrt{3}\right)\left(1+x-y/\sqrt{3}\right)\left|\sum_{i=2}^{m_3}c_{3,i}x^{i-2}\right|\\
    &\leq\max_{-1\leq x\leq 1}\frac{1}{4}(1-x^2)\left|\sum_{i=2}^{m_3}c_{3,i}x^{i-2}\right|
=\max_{(x,y)\in \hat{F}^3}|p|.
    \end{align*}
This proves the desired result.
\end{proof}

Finally, we can state the following quasi-optimal approximation property of $\Pi$.
\begin{theorem}
\label{local error theorem}
For every $m\in\mathbb{N}_+$, the interpolation operator $\Pi:C(\hat{T})\to \hat{P}$ satisfies
\begin{equation*}
\|f-\Pi f\|_{\infty,\hat{T}}\lesssim \cstab{m}\inf_{p\in \hat{P}}\|f-p\|_{\infty,\hat{T}}  \text{ with } \cstab{m}:=m\log m.
    \end{equation*}
\end{theorem}
\begin{proof}
The stability of $E$ and $\mathbf{E}$ in \eqref{eq:en} and \eqref{eq:ext} gives the stability of $\mathbf{E}E: C(\partial \hat{T})\to C(\hat{T})$:
\begin{equation}\label{eq:extension-stab}
    \|\mathbf{E}(E f)\|_{\infty,\hat{T}}\lesssim\|E f\|_{\infty,\partial \hat{T}}\lesssim\log m\|f\|_{\infty,\partial \hat{T}}.
\end{equation}
Together with the stability estimate of $I$ \eqref{eq:in}, it leads to
\begin{equation*}
    \begin{aligned}
        \|\Pi f\|_{\infty,\hat{T}}&\leq \|\mathbf{E}(E f)\|_{\infty,\hat{T}}+\|I(f-\mathbf{E}(Ef))\|_{\infty,\hat{T}} \\
        &\lesssim\log m\|f\|_{\infty,\partial \hat{T}}+m\|f-\mathbf{E}(Ef)\|_{\infty,\hat{T}} \\
        &\leq\log m\|f\|_{\infty,\partial \hat{T}}+m(\|f\|_{\infty,\hat{T}}
        +\|\mathbf{E}(E f)\|_{\infty,\hat{T}} ).
    \end{aligned}
\end{equation*}
In view of \eqref{eq:ext} and \eqref{eq:extension-stab}, we obtain
\begin{align*}
	\|\Pi f\|_{\infty,\hat{T}}&\lesssim
	\log m\|f\|_{\infty,\partial \hat{T}}+m\|f\|_{\infty,\hat{T}}
	+m\log m\|f\|_{\infty,\partial\hat{T}}
	\lesssim m\log m\|f\|_{\infty,\hat{T}}.
\end{align*}
Note that $\Pi p=p$ for all $p\in \hat{P}$, consequently, we arrive at
\begin{align*}
    \|f-\Pi f\|_{\infty,\hat{T}}&=\|(f-p)-\Pi(f-p)\|_{\infty,\hat{T}}
    \lesssim m\log m\|f-p\|_{\infty,\hat{T}}.
\end{align*}
Taking the infimum over all $p\in \hat{P}$ leads to the desired assertion.
\end{proof}
Furthermore, we can derive the best approximation in $\hat{P}$ in uniform norm.
\begin{lemma}[Best approximation in $\hat{P}$]\label{best approximation}
For all integer $0\leq j\leq \min (m_{1},m_{2},m_{3})$ and all $f\in W^{j+1,\infty}(\hat{T})$,
there exists a bounded constant $\mathrm{C}_{\rm bpa}$ independent of $j$ and $h_{\hat{T}}$ such that
\begin{equation}\label{ieq best approximation}
\inf_{p\in \hat{P}}\|f-p\|_{\infty,\hat{T}}\leq \mathrm{C}_{{\rm bpa}} \left(\frac{h_{\hat{T}}}{j+1}\right)^{j+1} |f|_{W^{j+1,\infty}(\hat{T})},
\end{equation}
where $|\cdot|_{W^{k,\infty}(D)}$ is the $W^{k,\infty}(D)$ semi-norm.
\end{lemma}
\begin{proof}
Recall the following well-known result \cite[Equation (5.4.16)]{spectral2006}: for all $k\geq 0$, there holds
\begin{align*}
\inf_{p\in P_k(\hat{T})}\|f-p\|_{\infty,\hat{T}}\leq \mathrm{C}_{\rm bpa} \left(\frac{h_{\hat{T}}}{k+1}\right)^{k+1} |f|_{W^{k+1,\infty}(\hat{T})}, \quad \forall f\in W^{k+1,\infty}(\hat{T}).
\end{align*}
By choosing $k:=\min (m_{1},m_{2},m_{3})$ and noting $P_k(\hat{T})\subset\hat{P}$, we establish the lemma.
\end{proof}
\begin{remark}[Boundedness of the constant $\mathrm{C}_{\rm bpa}$]\label{whitney-constant}
Note that the best approximation error $\inf_{p\in P_k(\hat{T})}\|f-p\|_{\infty,\hat{T}}$ can be bounded by the product of the so-called Whitney constant and the modulus of smoothness of order $k+1$ of a function $f\in C(\hat{T})$, and note also that the upper bound of the modulus of smoothness is given by $\left(\frac{h_{\hat{T}}}{k+1}\right)^{k+1}|f|_{W^{k+1,\infty}(\hat{T})}$ for $f\in W^{k+1,\infty}(\hat{T})$ \cite[Equation (2.4)]{johnen2006equivalence}. Furthermore, the Whitney constant is proved to be bounded \cite{kryakin1993exact, sendov1991theorem,dzyadyk2008theory, gilewicz2002boundedness, karaivanov2003nonlinear}. Consequently, the constant $\mathrm{C}_{\rm bpa}$ is bounded.
\end{remark}

\subsection{Global approximation error}
Now we analyze the convergence of Algorithm \ref{alg2}. The following two lemmas are standard.
\begin{lemma}[{\cite[Theorem 3.1.2]{ciarlet2002finite}}]\label{lemma 1}
Suppose $f\in W^{k,\infty}(T)$ for $k\in \mathbb{N}$ and some element $T\in\mathcal{T}_n$. Then $\hat{f}:=f\circ M_T$ belongs to $W^{k,\infty}(\hat{T})$, and there exists a constant $C(k)$ depending on $k$ such that
\begin{align*}
|\hat{f}|_{W^{k,\infty}(\hat{T})}&\leq C(k)\Vert A_T\Vert^k |f|_{W^{k,\infty}(T)},\\
|f|_{W^{k,\infty}(T)}&\leq C(k)\Vert A_T^{-1}\Vert^k |\hat{f}|_{W^{k,\infty}(\hat{T})},
\end{align*}
where $A_T$ is the matrix in the mapping $M_T$.
\end{lemma}
The next result gives bounds on $\Vert A_T\Vert$ and $\Vert A_T^{-1}\Vert$.
\begin{lemma}[{\cite[Theorem 3.1.3]{ciarlet2002finite}}]\label{lemma 2}
There holds $\Vert A_T\Vert \leq h_{T}\rho_{\hat{T}}^{-1}$ and $\Vert A_T^{-1}\Vert \leq h_{\hat{T}}\rho_{T}^{-1}$.
\end{lemma}

Next, we establish the element-wise approximation property of the global projection $\Pi_{\mathcal{T}_n}$.
\begin{theorem}[Element-wise approximation property of $\Pi_{\mathcal{T}_n}$]
\label{property}
Let $n\in\mathbb{N}$ be sufficiently large.
Let $\mathcal{T}_n$ with element-wise and facewise degrees $\{n_T\}_{T\in\mathcal{T}_n}$ and $\{m_F\}_{F\in\mathcal{F}_n}$ be generated by Algorithm \ref{alg1}, Definitions \ref{degree design strategy} and \ref{edge degree design strategy} and let the $hp$ approximation space be given by \eqref{def:vn}. Then for $j=1,\cdots,L-1$, there holds
\begin{align}
    \Vert \lambda_j-\Pi_{\mathcal{T}_n}\lambda_j\Vert_{\infty, T}&\lesssim   \mathrm{C}_{\rm bpa}h_T |\lambda_j|_{W^{1,\infty}(T)} , &&\text{for }T\in\mathcal{M}_n\label{error t2}\\
    \Vert \lambda_j-\Pi_{\mathcal{T}_n}\lambda_j\Vert_{\infty, T}&\lesssim \cstab{n_T} \mathrm{C}_{\rm bpa}C(\Tilde{n}) \left(\frac{\sqrt{3}h_T}{\Tilde{n}}\right)^{\Tilde{n}}|\lambda_j|_{W^{\tilde{n},\infty}(T)}, &&\text{for }T\in\mathcal{T}_n\backslash\mathcal{M}_n,\label{error t1}
\end{align}
where $\cstab{n_T}:=n_T\log n_T$, $C(\Tilde{n})$ is the constant in Lemma \ref{lemma 1} and $\Tilde{n}=\min(m_{F_T^1},m_{F_T^2},m_{F_T^3})+1$.
Furthermore, if the slope parameter $\mu$ satisfies the following condition
\begin{equation}\label{eq:mu-cond}
    \ln\left(\cstab{\mu \ell_T} {C(n_T+1)}C_{\lambda_j,T}^{{\left\lceil \mu \ell_T+1 \right\rceil}}\right)+{\left\lceil \mu \ell_T +1 \right\rceil} \left(\ln \sqrt{3}+{F(\ell_T)}-{\ln\left\lceil \mu \ell_T+1 \right\rceil}\right)
    \lesssim \ln\left(C_{\max}\right)+{F(1)}
\end{equation}
for all $T\in\mathcal{T}_n\backslash\mathcal{M}_n$, then
\begin{align*}
\max_{T\in \mathcal{T}_n\backslash\mathcal{M}_n}\Vert \lambda_j-\Pi_{\mathcal{T}_n}\lambda_j\Vert_{\infty, T}\lesssim
 \min_{T\in\mathcal{M}_n}\Vert\lambda_j-\Pi_{\mathcal{T}_n}\lambda_j\Vert_{\infty, T}.
\end{align*}
Here, $C_{\max}:=\max_{T\in\mathcal{M}_n}  |\lambda_j|_{W^{1,\infty}(T)}$, $C_{\lambda_j,T}^{\tilde{n}}:=|\lambda_j|_{W^{\tilde{n},\infty}(T)}$, and {$F(\ell):=\ln h_1-\frac{n-\ell}{2}\ln 2$} with $h_1$ being the initial mesh size.
\end{theorem}
\begin{proof}	
By the definition of the element-wise interpolation \eqref{eq:interp-global}, we obtain
\begin{align*}
\Vert \lambda_j-\Pi_{\mathcal{T}_n}\lambda_j\Vert_{\infty, T}
=\left\Vert\left(\lambda_j\circ M_T-\Pi(\lambda_j\circ M_T)\right)\circ M_T^{-1}\right\Vert_{\infty, T},\quad \forall T\in \mathcal{T}_n.
\end{align*}
Let $\hat{\lambda}_j:=\lambda_j\circ M_T$. Then we derive
\begin{align*}
\Vert \lambda_j-\Pi_{\mathcal{T}_n}\lambda_j\Vert_{\infty, T}
=\left\Vert\hat{\lambda}_j-\Pi\hat{\lambda}_j\right\Vert_{\infty, \hat{T}}.
\end{align*}
Together with Theorem \ref{local error theorem}, we derive
\begin{align}\label{eq:111}
\Vert \lambda_j-\Pi_{\mathcal{T}_n}\lambda_j\Vert_{\infty, T}
\lesssim \cstab{n_T}\inf_{p\in \hat{P}}\|\hat{\lambda}_j-p\|_{\infty,\hat{T}}.
\end{align}
Now the proof of the theorem consists of three steps.
First, if $T\in \mathcal{M}_n$, then $\lambda_j\in W^{1,\infty}(T)$ and $\hat{\lambda}_j\in W^{1,\infty}(\hat{T})$. Moreover, by Definition \eqref{degree design strategy}, we have $n_T=2$. Thus, $\cstab{n_T}\leq 3$.  By Lemma \ref{best approximation}, we derive
\begin{equation*}
    \inf_{p\in \hat{P}}\|\hat{\lambda}_j-p\|_{\infty,\hat{T}}\leq \mathrm{C}_{\rm bpa} h_{\hat{T}} |\hat{\lambda}_j|_{W^{1,\infty}(\hat{T})}.
\end{equation*}
Then combining with Lemmas \ref{lemma 1} and \ref{lemma 2}, this further leads to
\begin{align*}
\Vert \lambda_j-\Pi_{\mathcal{T}_n}\lambda_j\Vert_{\infty, T}
&\lesssim\mathrm{C}_{\rm bpa} h_{\hat{T}} |\hat{\lambda}_j|_{W^{1,\infty}(\hat{T})}
\lesssim\mathrm{C}_{\rm bpa}h_T|{\lambda_j}|_{W^{1,\infty}({T})}.
\end{align*}
This proves the assertion \eqref{error t2}. Second, if $T\in \mathcal{T}_n\backslash\mathcal{M}_n$, then $\hat{\lambda}_j:=\lambda_j\circ M_T\in W^{\Tilde{n},\infty}(T)$, and \eqref{eq:111} and Lemma \ref{best approximation} yield
\begin{equation*}
    \inf_{p\in \hat{P}}\|\hat{\lambda}_j-p\|_{\infty,\hat{T}}\leq \mathrm{C}_{\rm bpa} (h_{\hat{T}}/{\Tilde{n}})^{\Tilde{n}} |\hat{\lambda}_j|_{W^{\Tilde{n},\infty}(\hat{T})}.
\end{equation*}
Together with Lemmas \ref{lemma 1} and \ref{lemma 2} and noting $\frac{h_{\hat{T}}}{\rho_{\hat{T}}}=\sqrt{3}$, we obtain
\begin{align*}
\Vert \lambda_j-\Pi_{\mathcal{T}_n}\lambda_j\Vert_{\infty, T}
&\lesssim\cstab{n_T}\mathrm{C}_{\rm bpa} ({h_{\hat{T}}}/{\Tilde{n}})^{\Tilde{n}} |\hat{\lambda}_j|_{W^{\Tilde{n},\infty}(\hat{T})}\\
&\leq \cstab{n_T}\mathrm{C}_{\rm bpa}C(\Tilde{n}) ({h_{\hat{T}}}/{\Tilde{n}})^{\Tilde{n}}(h_T\rho_{\hat{T}}^{-1})^{\Tilde{n}}|\lambda_j|_{W^{\tilde{n},\infty}(T)}\\
&=\cstab{n_T}\mathrm{C}_{\rm bpa}C(\Tilde{n})({\sqrt{3}h_T}/{\Tilde{n}})^{\Tilde{n}}|\lambda_j|_{W^{\tilde{n},\infty}(T)}.
\end{align*}
This proves the assertion \eqref{error t1}.
Finally, we prove the last statement that if the slope parameter $\mu$ satisfies \eqref{eq:mu-cond}, then for all elements $T\in\mathcal{T}_n^R$ with  $n_T=\left\lceil \mu \ell_T \right\rceil$, the right-hand side of \eqref{error t1} is no greater than the right-hand side of \eqref{error t2} for some element in $\mathcal{M}_n$.
Definition \ref{edge degree design strategy} implies
$n_T\geq m_{F_T^i}$, for $i=1,2,3$,
for all $T\in \mathcal{T}_n\backslash\mathcal{M}_n$. Moreover, if $m_{F_T^i}\leq n_T$, there exists $\Tilde{T}\in\mathcal{F}_F$ such that $m_{F_T^i}=n_{\Tilde{T}}$. This motivates setting {$\Tilde{n}=n_T+1$} in \eqref{error t1} when analyzing the choice of $\mu$.
Algorithm \ref{alg1} and Definition \ref{layer} indicate the diameters of elements satisfy {$ h_T\sim h_12^{-(n-\ell_T)/2}$}.
Thus, ignoring all the absolute constants, the right hand side of \eqref{error t1} becomes
\begin{equation}\label{eq: right hand 1}
    \cstab{n_T}\mathrm{C}_{\rm bpa}{C(n_T+1)}({\sqrt{3}h_12^{-(n-\ell_T)/2}}/{{(n_T+1)}})^{{n_T+1}}|\lambda_j|_{{W^{n_T+1,\infty}(T)}},
\end{equation}
and the right-hand side of \eqref{error t2} is
\begin{equation}\label{eq: right hand 2}
    \mathrm{C}_{\rm bpa}h_12^{-(n-1)/2} |\lambda_j|_{W^{1,\infty}(T)}.
\end{equation}
Then the condition \eqref{eq:mu-cond} implies the maximum value of \eqref{eq: right hand 2} for all $T\in\mathcal{M}_n$ is bounded above by the minimum value of \eqref{eq: right hand 2} for all $T\in\mathcal{T}_n\backslash\mathcal{M}_n$, and this completes the proof.
\end{proof}

\begin{remark}\label{degree choice}
Note that the constant $\mathrm{C}_{\rm bpa}$ is bounded, cf Remark \ref{whitney-constant}. Assume that $h_T$ and $h_T^{\Tilde{n}}$ dominate the right-hand sides of \eqref{error t2} and \eqref{error t1} respectively, with a small enough initial mesh size, $\mu=1$ or $\mu=0.5$ is sufficient for \eqref{eq:mu-cond}. Then the point with the largest error lies in $\mathcal{M}_n$, as is numerically verified.
\end{remark}

Next we discuss the approximation error under the assumption that the maximum error exists in some element belonging to $\mathcal{M}_n$. Since the mesh $\mathcal{T}_n$ generated by the procedure is geometrically refined towards branch points with $n$ layers, we derive a convergence rate according to the classification of branch points.

\begin{figure}[hbt!]
\centering
\subfigure[One singular point $\frac{1}{a}(\frac{\pi}{3},\frac{\pi}{3})$]{\label{refine demo finite singular}
\includegraphics[width=3.7cm,trim={4cm 0.2cm 4cm 1cm},clip]{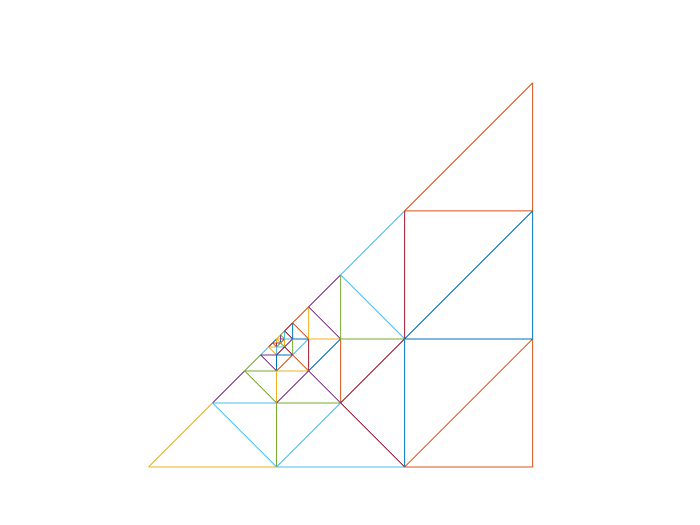}
}%
\subfigure[Two singular points $\frac{1}{a}(\frac{\pi}{3},\frac{\pi}{3})$ and $\frac{1}{a}(\frac{2\pi}{3},\frac{\pi}{3})$]{\label{refine demo finite singular2}
\includegraphics[width=3.7cm,trim={4cm 0.2cm 4cm 1cm},clip]{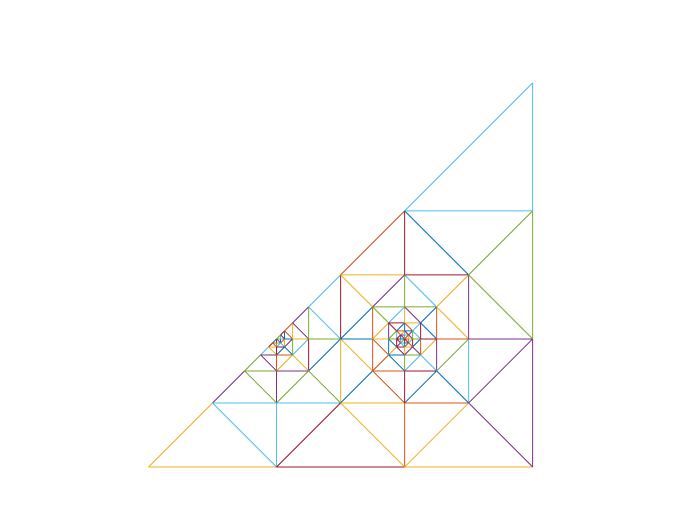}
}%
\subfigure[One singular line $k_2=k_1$]{\label{refine demo infinite singular}
\includegraphics[width=3.7cm,trim={4cm 0.2cm 4cm 1cm},clip]{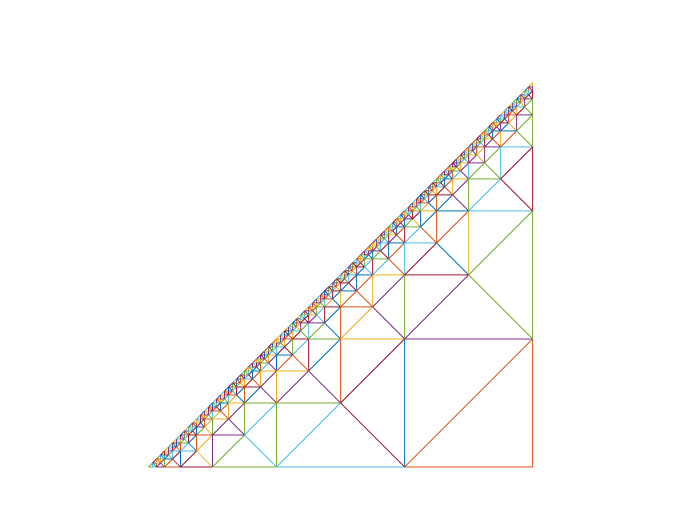}
}%
\subfigure[Two singular lines $k_2=k_1$ and $k_2+k_1+\frac{\pi}{a}=0$ with $\frac{2}{a}\leq k_1\leq \frac{3}{a}$]{\label{refine demo infinite singular2}
\includegraphics[width=3.7cm,trim={4cm 0.2cm 4cm 1cm},clip]{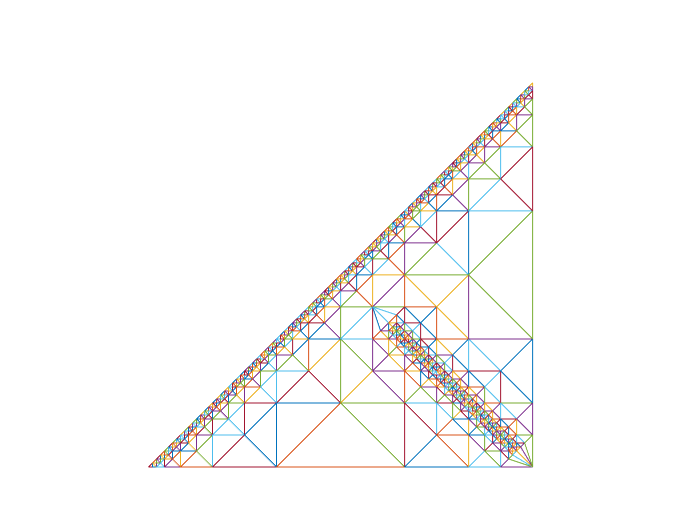}
}%
\caption{Algorithm \ref{alg1} with finite (the first two) and infinite (the last two) singularities.}\label{refine demo}
\end{figure}

First, if the band functions we consider have a finite number of branch points, then the number of elements $\# \mathcal{T}_{n}=\mathcal{O}(n)$ when $n$ is sufficiently large, leading to the number of sampling points $N=\mathcal{O}(n^3)$. See Figure \ref{refine demo} for an illustration, in which the left two figures are the mesh refinement in case of a finite number of branch points in the IBZ of the square lattice. In this case, we derive an exponential convergence rate,

\begin{theorem}[Convergence rate for $hp$ element-wise interpolation with finite branch points]\label{exp}
Let $n\in\mathbb{N}$ be sufficiently large.
Let $\mathcal{T}_n$ with element-wise and facewise degrees $\{n_T\}_{T\in\mathcal{T}_n}$ and $\{m_F\}_{F\in\mathcal{F}_n}$ be generated by Algorithm \ref{alg1}, Definitions \ref{degree design strategy} and \ref{edge degree design strategy}, the $hp$ approximation space be defined in \eqref{def:vn}, and the slope parameter $\mu$ satisfy \eqref{eq:mu-cond}.
If the band functions we consider have a finite number of branch points, then there holds
\begin{equation*}
\Vert \lambda_j-\Pi_{\mathcal{T}_n}\lambda_j\Vert_{\infty,\mathcal{B}}\lesssim \exp(-b N^{\frac{1}{3}})| \lambda_j|_{W^{1,\infty}(\mathcal{B})}\text{ for } j=1,\cdots,L-1.
\end{equation*}
Here, the constant $b>0$ is independent of the number of sampling points $N$.
\end{theorem}
\begin{proof}
Theorem \ref{refining mechanism} indicates that the set of singular points $\mathbf{S}$ as a subset
of the branch points is marked at the $n$th iteration, i.e., $\mathbf{S}\subset\mathcal{M}_n$.
Algorithm \ref{alg1} implies that the diameter of each element $T\in\mathcal{M}_n$ satisfies
$h_T\sim 2^{-n/2}$ for all $T\in\mathcal{M}_n$. The choice of $\mu$ (i.e., the minimum
value satisfying \eqref{eq:mu-cond}) ensures that the element-wise interpolation error measured
in the uniform norm attains its maximum over $\mathcal{M}_n$. This and Theorem \ref{property} imply
\begin{equation}\label{error 1}
\Vert \lambda_j-\Pi_{\mathcal{T}_n}\lambda_j\Vert_{\infty,\mathcal{B}}\lesssim 2^{-n/2}| \lambda_j|_{W^{1,\infty}(\mathcal{B})}\text{ for } j=1,\cdots,L-1.
\end{equation}
Since the number of sampling points $N=\mathcal{O}(n^3)$ for sufficiently large $n$, the desired assertion follows.
\end{proof}

Next we discuss the convergence rate when the number of branch points is infinite; see Figure \ref{refine demo} for an illustration of the mesh refinement. Asymptotically, we have the number of elements $\# \mathcal{T}_{n}=\mathcal{O}(2^{\frac{n}{2}})$ and the number of sampling points $N\sim 2^{\frac{n}{2}}$. Similarly, we derive the following convergence rate for $hp$ element-wise interpolation with infinite branch points.
\begin{theorem}\label{algebraic}
Under the conditions of Theorem \ref{exp}, if the number of branch points is infinite, then
\begin{equation*}
\Vert \lambda_j-\Pi_{\mathcal{T}_n}\lambda_j\Vert_{\infty,\mathcal{B}}\lesssim N^{-1}| \lambda_j|_{W^{1,\infty}(\mathcal{B})}\text{ for } j=1,\cdots,L-1.
\end{equation*}
\end{theorem}

\begin{remark}
With $n\in\mathbb{N}$ sufficiently large, the approximation error of the band function is given by 
\begin{align*}
\Vert \omega_j-c\cdot\sqrt{\Pi_{\mathcal{T}_n}\lambda_j}\Vert_{\infty,\mathcal{B}}
&={c^2}\cdot\left\|\frac{ \lambda_j-\Pi_{\mathcal{T}_n}\lambda_j}{\omega_j+c\cdot\sqrt{\Pi_{\mathcal{T}_n}\lambda_j}}\right\|_{\infty,\mathcal{B}}
\\
&\leq{c\cdot}\frac{\Vert \lambda_j-\Pi_{\mathcal{T}_n}\lambda_j\Vert_{\infty,\mathcal{B}}}{2\sqrt{\lambda_j(\tilde{\mathbf{k}})}-\sqrt{\text{err}}}, 
\end{align*}
where $\tilde{\mathbf{k}}=\arg\max_{{\mathbf{k}\in\mathcal{B}}}{| \omega_j(\mathbf{k})-c\cdot\sqrt{\Pi_{\mathcal{T}_n}\lambda_j(\mathbf{k})}|}$ and $\text{err}=\Vert \lambda_j-\Pi_{\mathcal{T}_n}\lambda_j\Vert_{\infty,\mathcal{B}}$. 

Thus, the decay rates of {$\Vert \omega_j-c\cdot\sqrt{\Pi_{\mathcal{T}_n}\lambda_j}\Vert_{\infty,\mathcal{B}}$ and $\Vert \lambda_j-\Pi_{\mathcal{T}_n}\lambda_j\Vert_{\infty,\mathcal{B}}$} are roughly of the same order.
\end{remark}

\begin{remark}
Theorems \ref{exp} and \ref{algebraic} hold asymptotically when the number of iterations $n$ is sufficiently large. In practice, the largest allowed iteration $n_{\max}$ of Algorithm \ref{alg1} is not large. Then the number $N$ of sampling points is $N=\mathcal{O}( 2^{\frac{n}{\theta}})$ with $\theta>2$, which leads to a faster convergence rate. Indeed, we observe a second-order convergence rate in numerical experiments.
\end{remark}

Last, we discuss the tolerance $\operatorname{tol}_2$ and stopping criterion $n_{\max}$ in Algorithm \ref{alg1}. By \eqref{error 1}, we have
 \begin{equation*}
 {\Vert \lambda_j-\Pi_{\mathcal{T}_n}\lambda_j\Vert_{\infty,\mathcal{B}}}\lesssim2^{-n/2}| \lambda_j|_{W^{1,\infty}(\mathcal{B})}\text{ for } j=1,\cdots,L-1.
\end{equation*}
Thus, $\operatorname{tol}_2$ and $n_{\max}$ should satisfy $2^{-n_{\max}/2}\simeq
\operatorname{tol}_2\simeq \epsilon$, i.e., comparable with the target precision $\epsilon$.

\section{Numerical experiments}\label{sec:simulation}
\begin{figure}[hbt!]
\centering
\subfigure[Square lattice \#DOFs=2107]{\label{mesh1}
\includegraphics[width=0.28\textwidth,trim={2cm 1cm 2cm 2cm},clip]{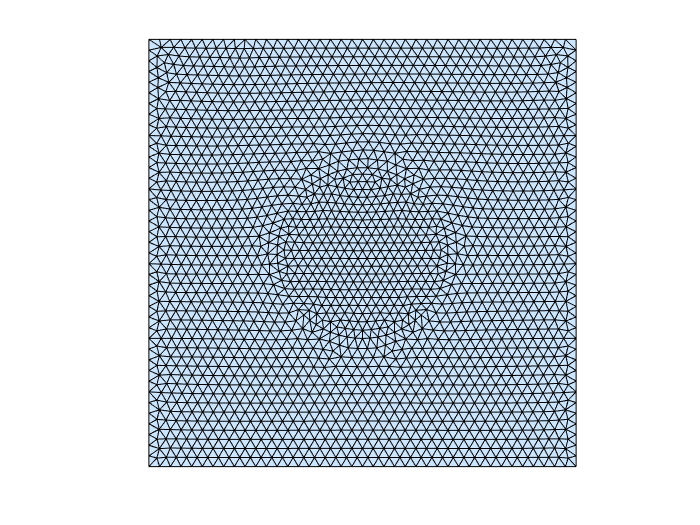}}\qquad
\subfigure[Hexagonal lattice \#DOFs=1781]{\label{mesh2}\includegraphics[width=0.3\textwidth,trim={1cm 1.7cm 1cm 1.3cm},clip]{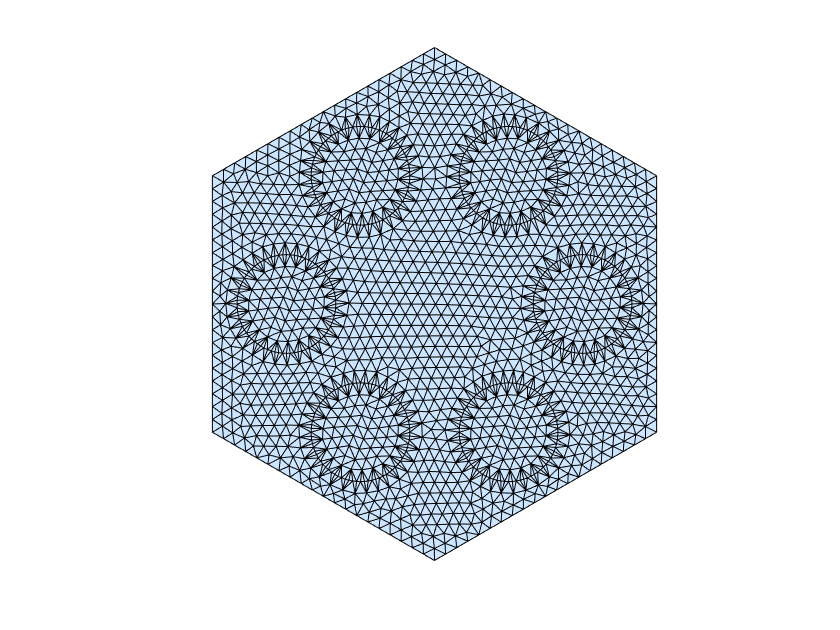}}
\caption{Discretization of the unit cell $\Omega$ by \textit{distmesh2d}.}
\end{figure}

To show the performance of Algorithm \ref{alg2}, we consider 2D PhCs with a square unit cell in
Figure \ref{lattice} and a hexagonal unit cell in Figure \ref{lattice3} as in \cite{wang2023dispersion}.
Given  $\mathbf{k}\in\mathcal{B}$, we solve the Helmholtz eigenvalue problem
\eqref{simply} using the conforming Galerkin FEM. Due to the high contrast between the relative
permittivity of the circular medium and the external one, we employ a fitted mesh generated by
\textit{distmesh2d} \cite{persson2004simple}; 
see Figures \ref{mesh1} and \ref{mesh2} for the mesh $\mathcal{T}_h$ with
a mesh size $h=0.025a$ for square lattice and $h=0.05a$ for hexagonal lattice. The associated
conforming piecewise affine space is given by
\begin{align*}
V_h=\{v_h\in C(\overline{\Omega}):v_h|_K\in P_1(K)\,\,\forall K\in\mathcal{T}_h\,\}.
\end{align*}
Given $\mathbf{k}\in\mathcal{B}$, the Galerkin FEM approximation to problem \eqref{simply} is to
find non-trivial eigenpair $(\lambda_h,u_h) \in (\mathbb{R},V_h)$ such that
\begin{equation}\label{fem variational}
    \left\{\begin{aligned}
&a(u_h,v_h)=\lambda_h b(u_h,v_h) , \quad \forall v \in V_h,\\
&\|u_h\|_{L^2_{\beta}(\Omega)}=1.
\end{aligned}\right.
\end{equation}

\begin{figure}[hbt!]
\centering
\subfigure[Square lattice TE mode]{\label{square close points}
\includegraphics[width=6cm,trim={10cm 1.2cm 8cm 1.0cm},clip]{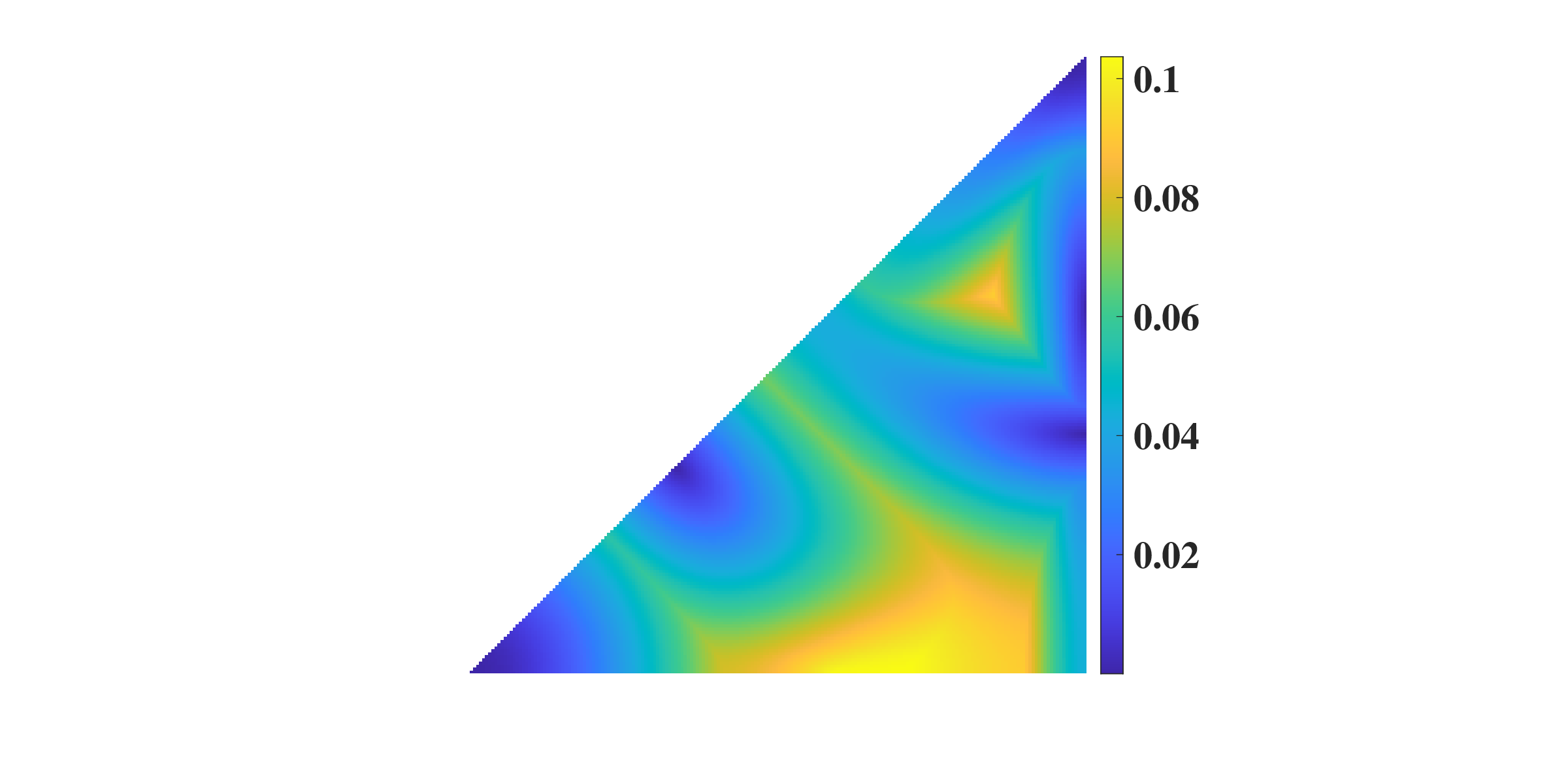}
}%
\quad
\subfigure[Hexagonal lattice TE mode]{\label{hex close points}
\includegraphics[width=7.5cm,trim={3cm 1.2cm 2cm 1.0cm},clip]{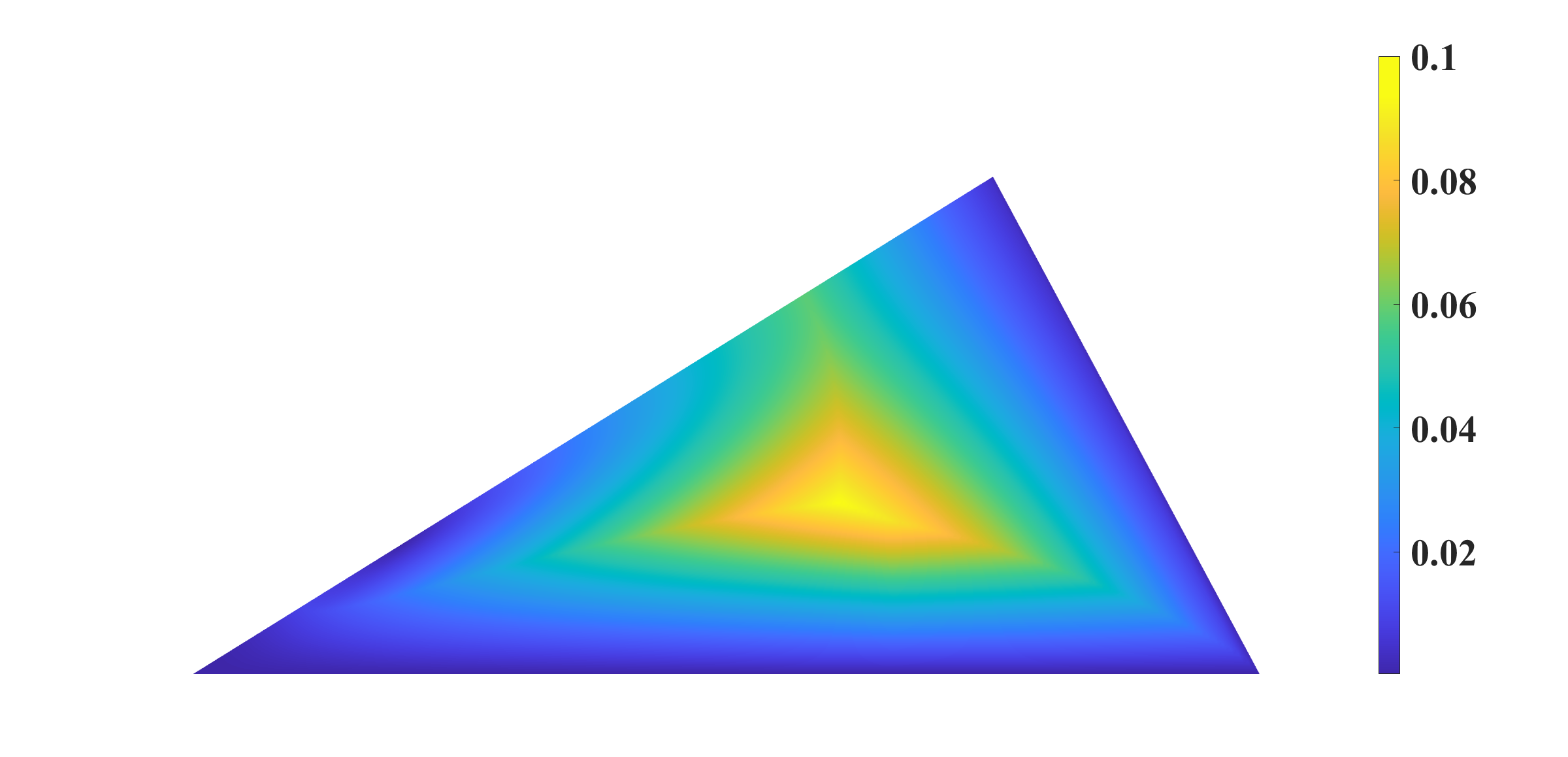}
}
\caption{Locations with distance between two adjacent band functions below $0.01$ are marked blue.}\label{close points}
\end{figure}

In both cases, IBZ is a triangle. To obtain the reference solution with sufficient accuracy,
we discretize $\mathcal{B}$ with 20503 evenly distributed points (denoted by $\hat{\mathcal{B}}$)
to compute the interpolation error and the accuracy of the proposed method.
The pointwise relative error is then computed on $\hat{\mathcal{B}}$ by
 \begin{align*}
    e_i(\mathbf{k}):=\frac{\vert \omega_i(\mathbf{k}) -\hat \omega_i(\mathbf{k})\vert}{{\omega_i}(\mathbf{k})} \quad \text{ for }\mathbf{k}\in\hat{\mathcal{B}} \text{ and } i=1,\cdots, 5.
 \end{align*}
Here, $\omega_i(\mathbf{k})$ is the $i$th reference band function obtained directly by the Galerkin FEM over $\hat{\mathcal{B}}$ using the same mesh on the unit cell $\Omega$, and $\hat\omega_i(\mathbf{k})$ is by the element-wise polynomial interpolation. Since the first few low-frequency band functions are typically of practical interest \cite{joannopoulos2008molding}, we show the performance of the numerical scheme using only the first five band functions. We measure the accuracy of the methods using maximum relative error
\begin{align*}
\operatorname{error}_{\infty}:=\max_{i=1,\cdots,5}\max_{\mathbf{k}\in\hat{\mathcal{B}}}\left|e_i(\mathbf{k})\right|.
\end{align*}
Figure \ref{close points} shows the location of singularities in the band functions with the distance between adjacent band functions. The areas where adjacent band functions are close or coincident are indicated in dark blue. These regions are targeted by Algorithm \ref{alg1}, allowing to detect the locations of singularities.

\subsection{Numerical tests with square lattice in Figure \ref{lattice}}
We first test Algorithm \ref{alg2} on the square lattice with the TE mode and compare the performance
of GPI with Fekete points and Chebyshev points of the second kind (Cheb2).

\begin{figure}[hbt!]
\centering
\subfigure[$\mu=1$]{\label{loglog_six2}
\includegraphics[width=8cm,trim={8cm 0.0cm 8cm 0.8cm},clip]{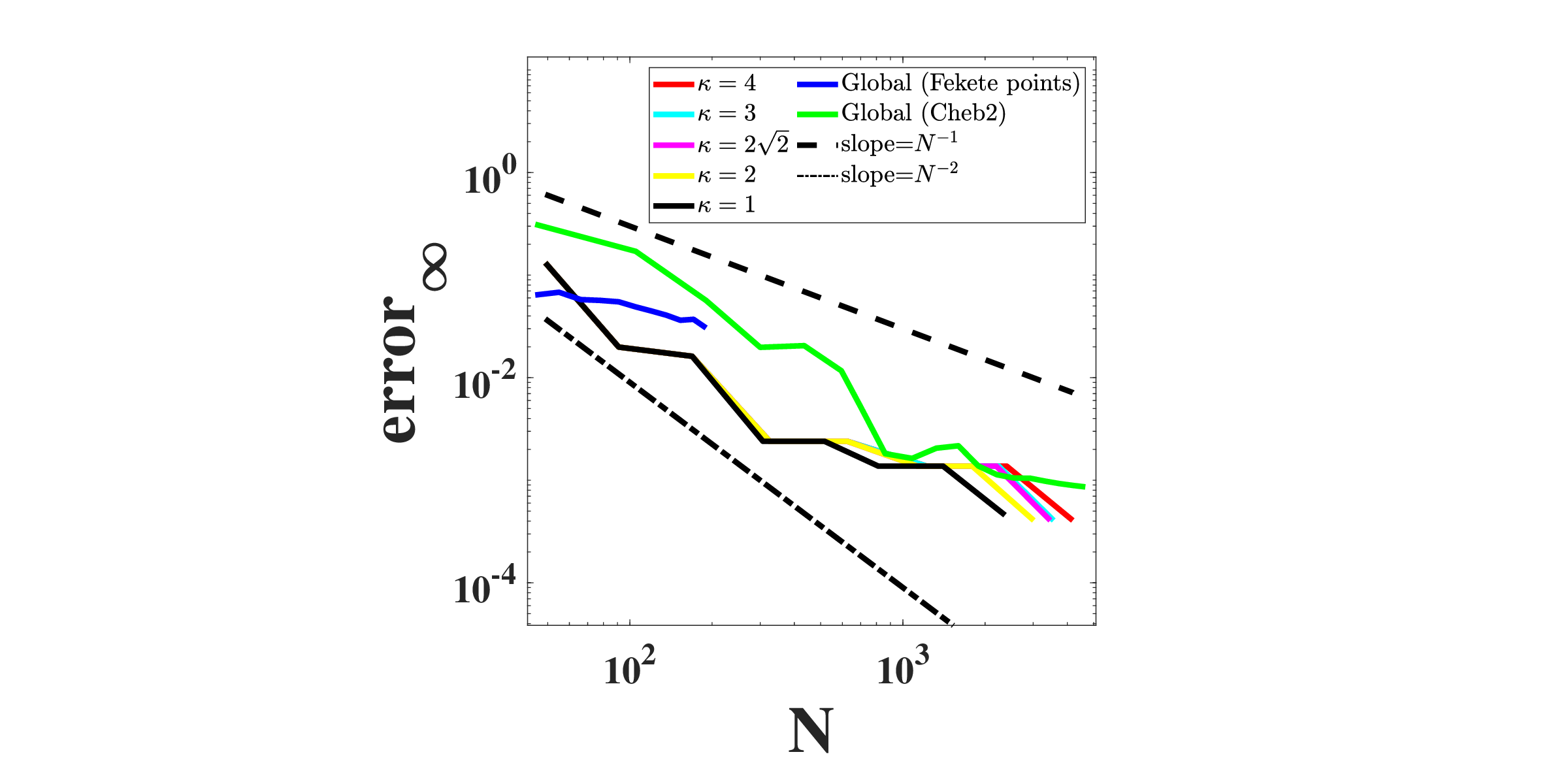}
}%
\quad
\subfigure[$\mu=0.5$]{\label{loglog_six2_mu2}
\includegraphics[width=8cm,trim={8cm 0.0cm 8cm 0.8cm},clip]{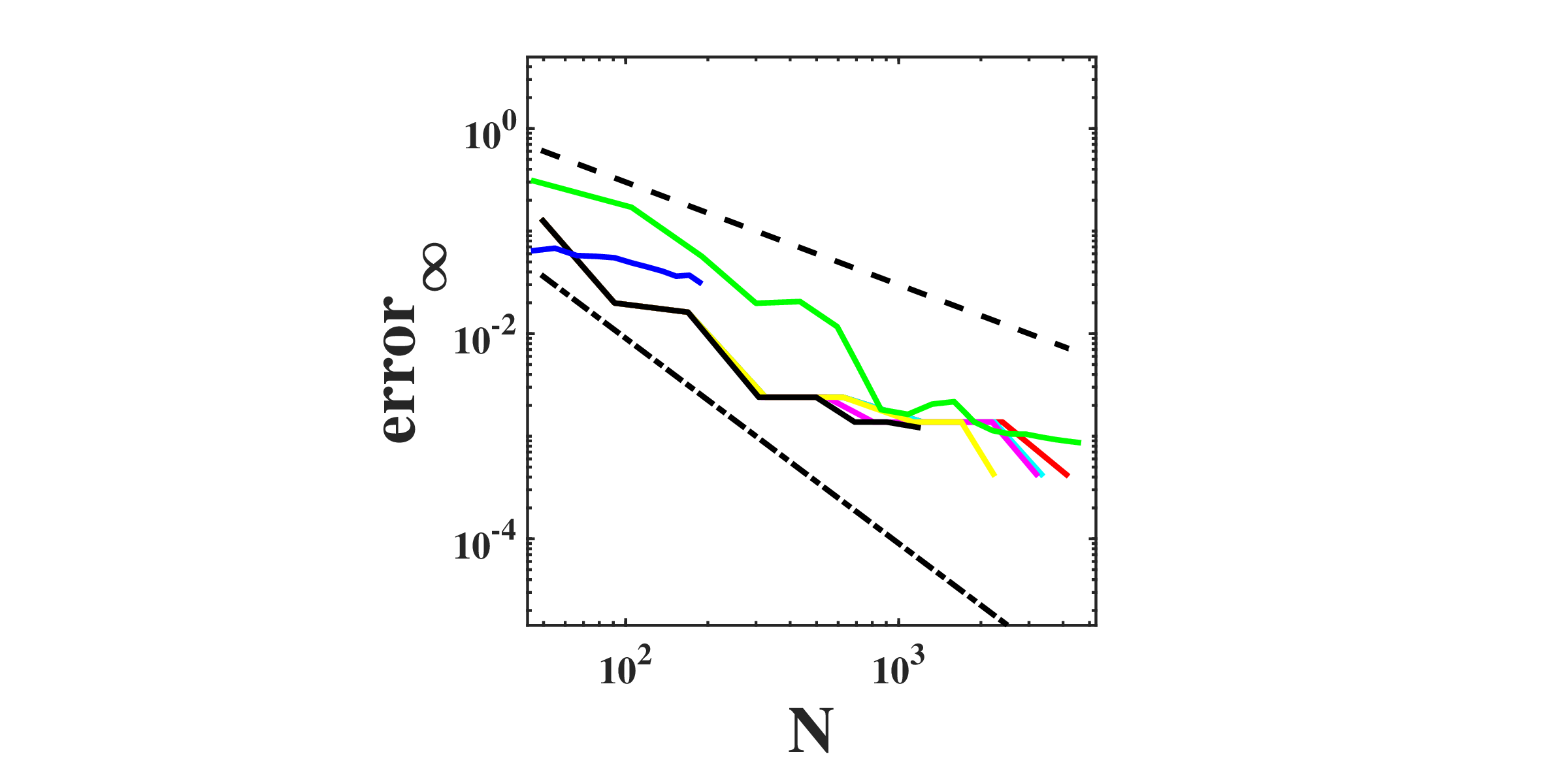}
}%
\quad
\centering
\caption{Maximum relative error with respect to the number of sampling points.}\label{loglog_six2 X}
\end{figure}

We present in Figure \ref{loglog_six2 X} the decay of
$\operatorname{error}_{\infty}$ against $N$ on a log-log scale of Algorithm \ref{alg2}
with $\kappa\in\{1,2,2\sqrt{2},3,4\}$ and $\mu\in\{1,0.5\}$ as suggested in Remark
\ref{degree choice}. The performance of GPI with Fekete points and Cheb2 are also provided
to facilitate comparison. Since Fekete points for polynomial spaces with (total) degrees
higher than 18 are still unknown \cite{taylor2000algorithm}, GPI with Fekete points is
restricted to $N\leq 190$. We observe that Algorithm \ref{alg2} outperforms GPI for all
cases. Furthermore, an algebraic convergence rate between first order and second order
is observed for any parameters of $\kappa$ and $\mu$.


\begin{figure}[hbt!]
\centering
\subfigure[loop 2]{\label{V5loop2}
\includegraphics[width=3.7cm,trim={2cm 1.5cm 2cm 2cm},clip]{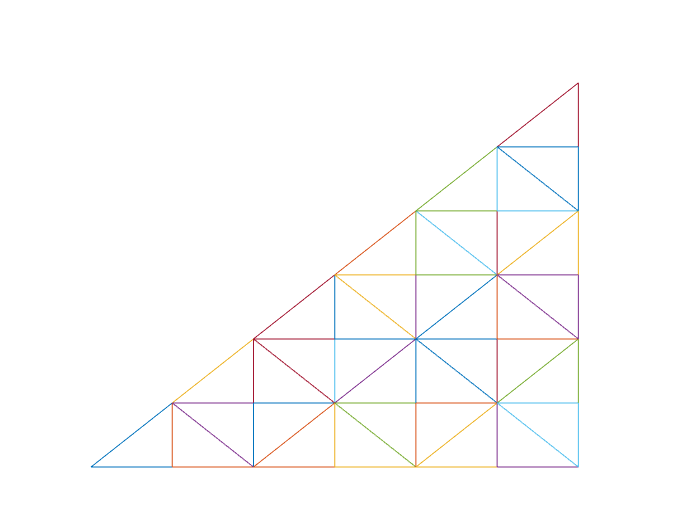}
}%
\quad
\subfigure[loop 4]{\label{V5loop4}
\includegraphics[width=3.7cm,trim={2cm 1.5cm 2cm 2cm},clip]{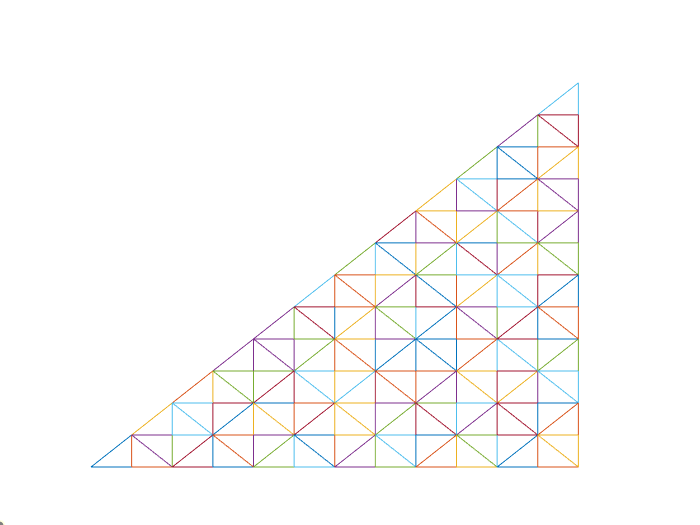}
}%
\quad
\subfigure[loop 6]{\label{V5loop6}
\includegraphics[width=3.7cm,trim={2cm 1.5cm 2cm 2cm},clip]{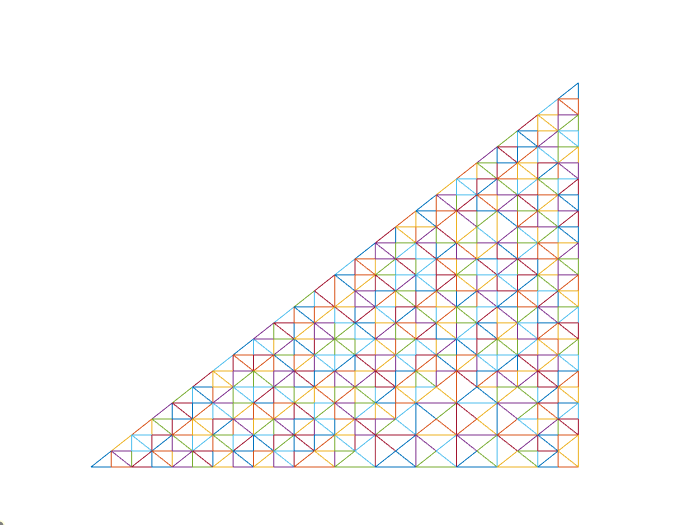}
}%
\quad
\subfigure[loop 8]{\label{V5loop8}
\includegraphics[width=3.7cm,trim={2cm 1.5cm 2cm 2cm},clip]{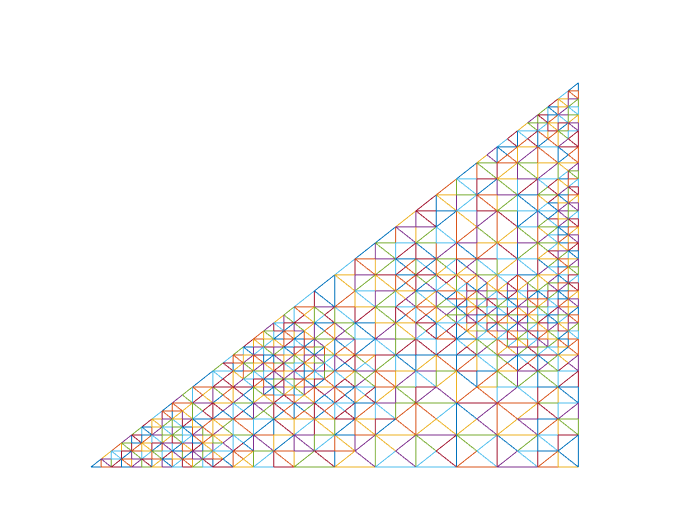}
}%
\quad
\subfigure[loop 2]{\label{V2loop2}
\includegraphics[width=3.7cm,trim={2cm 1.5cm 2cm 2cm},clip]{figs_NEW/V5loop2}
}%
\quad
\subfigure[loop 4]{\label{V2loop4}
\includegraphics[width=3.7cm,trim={2cm 1.5cm 2cm 2cm},clip]{figs_NEW/V5loop4}
}%
\quad
\subfigure[loop 6]{\label{V2loop6}
\includegraphics[width=3.7cm,trim={2cm 1.5cm 2cm 2cm},clip]{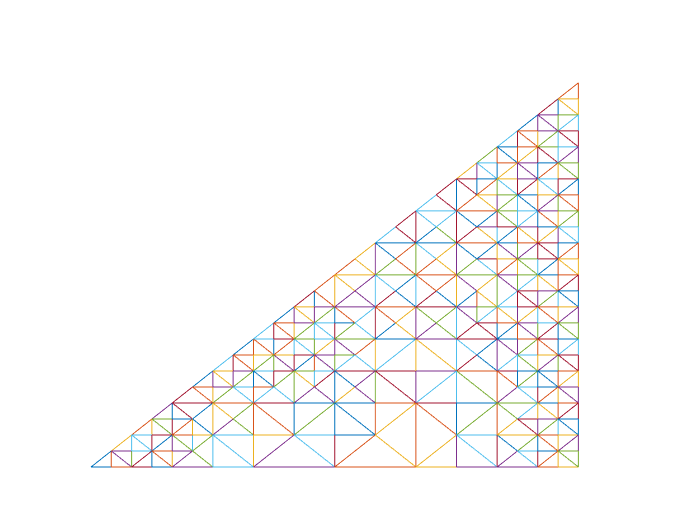}
}%
\subfigure[loop 8]{\label{V2loop8}
\includegraphics[width=3.7cm,trim={2cm 1.4cm 2cm 2cm},clip]{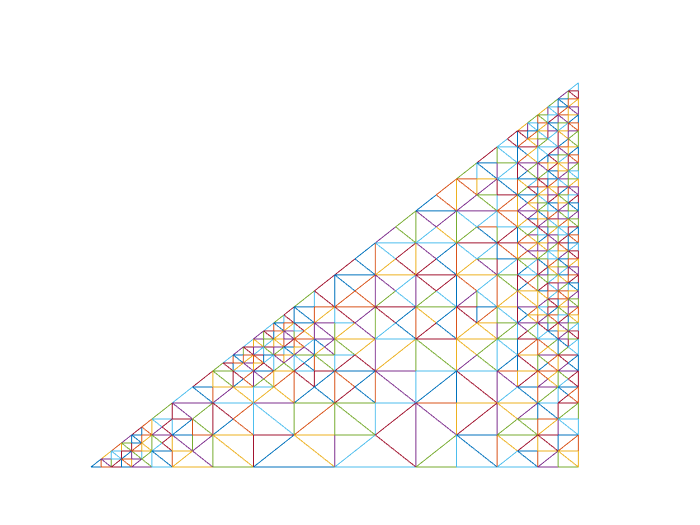}
}%
\centering
\caption{Algorithm \ref{alg1} with $\kappa=2\sqrt{2}$ (the top row) and $\kappa=2$ (the bottom row).}\label{VXloopX}
\end{figure}
Figure \ref{VXloopX} displays the adaptive meshes generated by Algorithm \ref{alg1} with
$\kappa=2\sqrt{2}$ and $\kappa=2$, respectively. Comparing with Figure \ref{square close points},
one observes that Algorithm \ref{alg1} adaptively refines the elements with crossings. Note
that larger $\kappa$ indicates more elements being marked and refined in each loop. The
initial mesh has little effect on the resulting adaptive mesh since meshes overlap in
the first several loops.

\begin{figure}[hbt!]
\centering
\subfigure[loop 8]{\label{kappa2q2 error 1}
\includegraphics[width=6cm,trim={1.1cm 0.8cm 0cm 0cm},clip]{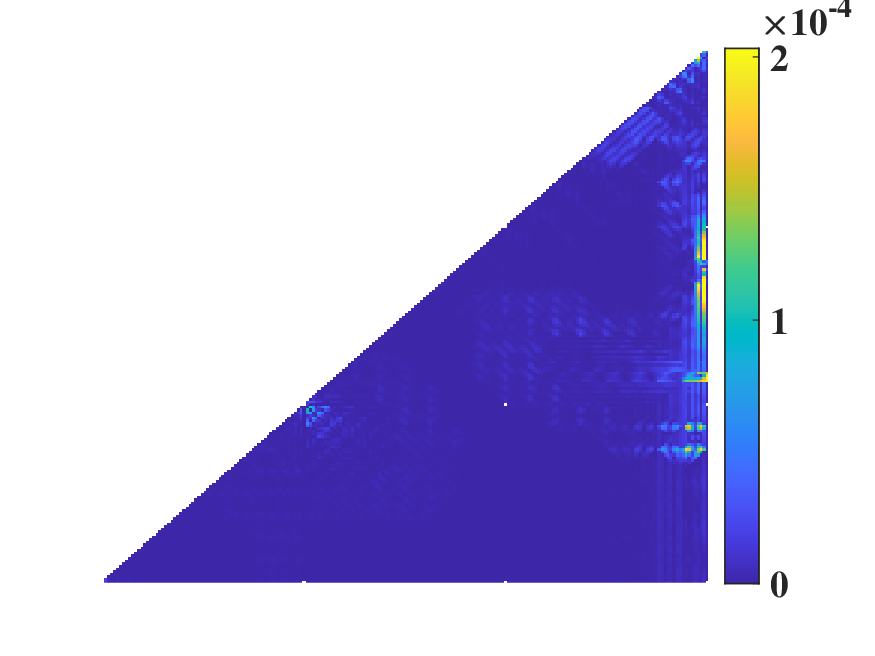}
}%
\quad
\subfigure[loop 8]{\label{kappa2q2 error 2}
\includegraphics[width=6cm,trim={1.1cm 0.8cm 0cm 0cm},clip]{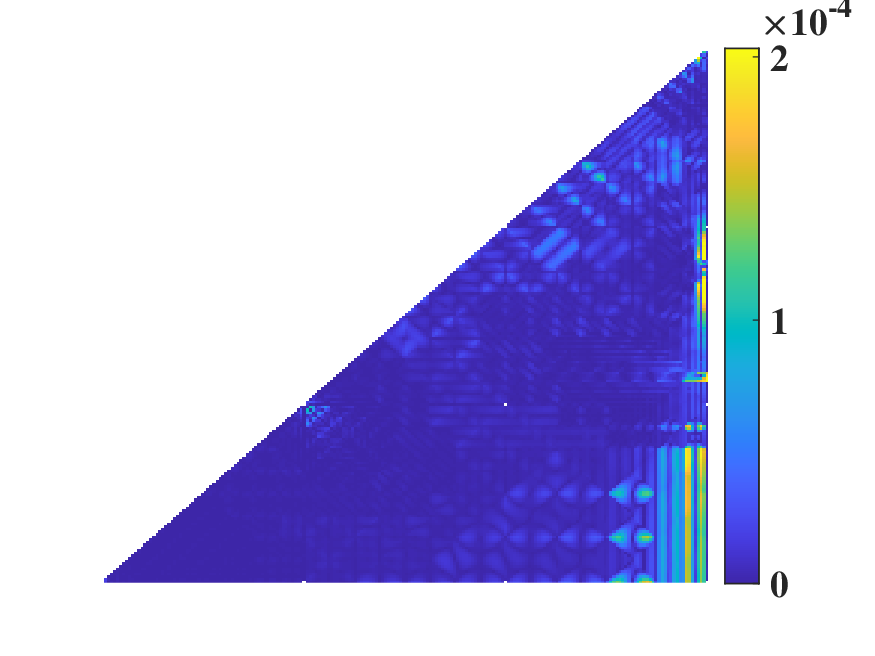}
}%
\quad
\centering
\caption{$\kappa=2\sqrt{2}$: pointwise relative error $e_i(\mathbf{k})$ under $\mu=1$ (left),  $\mu=0.5$ (right).}\label{kappa2q2 error}
\end{figure}

\begin{figure}[hbt!]
\centering
\subfigure[loop 8]{\label{kappa2 error 1}
\includegraphics[width=6cm,trim={1.1cm 0.8cm 0cm 0cm},clip]{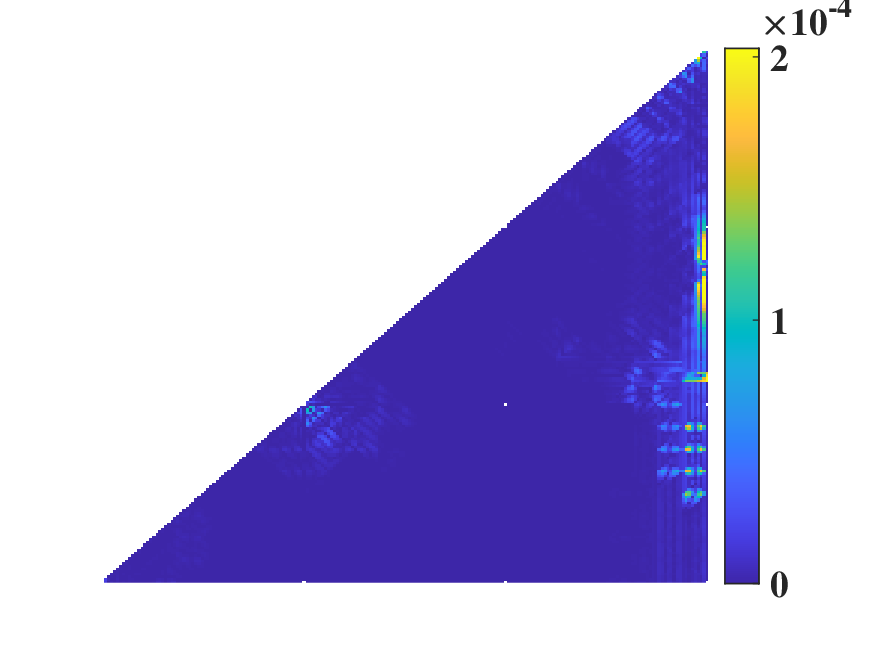}
}%
\quad
\subfigure[loop 8]{\label{kappa2 error 2}
\includegraphics[width=6cm,trim={1.1cm 0.8cm 0cm 0cm},clip]{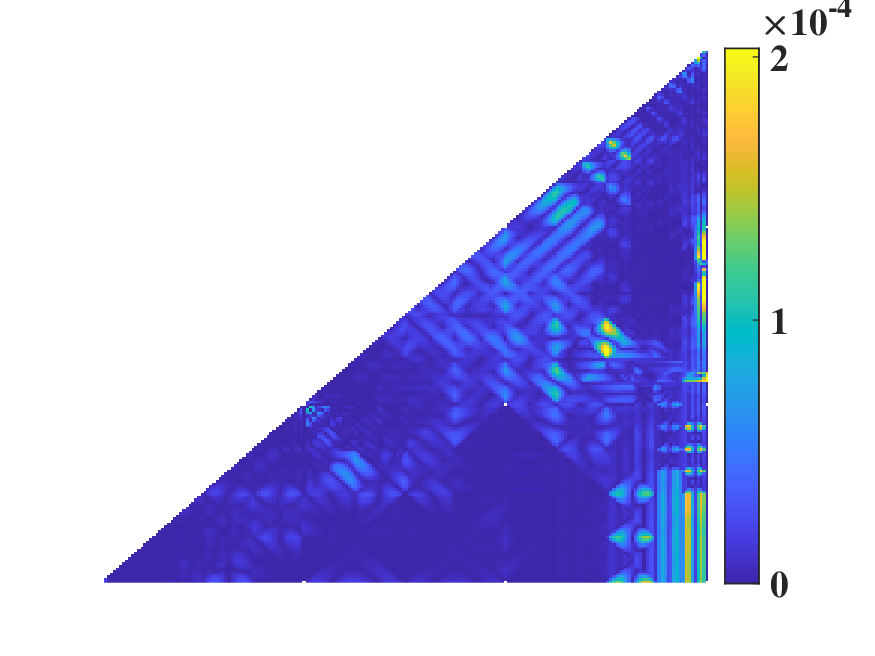}
}%
\quad
\centering
\caption{$\kappa=2$: pointwise relative error $e_i(\mathbf{k})$ under $\mu=1$ (left),  $\mu=0.5$ (right).}\label{kappa2 error}
\end{figure}

Figures \ref{kappa2q2 error} and \ref{kappa2 error} display the pointwise relative error
$e_i(\mathbf{k})$ of approximating the first five band functions under $\kappa=2\sqrt{2}$
and $\kappa=2$ and $\mu=1$ and $\mu=0.5$. Due to a lack of information about the seventh
band function, we are unable to determine the positions of all the singular points of the
sixth band function. Therefore, the error analysis only takes into account the first five
band functions, even though we interpolate the first six band functions. From Figures
\ref{kappa2q2 error 1} and \ref{kappa2 error 1}, one observes that with $\mu=1$ (i.e., the
total polynomial degree on each element is consistent with its layer), the pointwise relative
error $e_i(\mathbf{k})$ dominates in the singular element patch $\mathcal{T}_n^{S}$ with
$n=8$. This confirms Theorem \ref{refining mechanism}, and \eqref{eq:mu-cond}. However, for
a smaller slope parameter $\mu=0.5$, the pointwise relative
error $e_i(\mathbf{k})$ in Figures \ref{kappa2q2 error 2} and \ref{kappa2 error 2} can take
large values in the regular element patch $\mathcal{T}_n^{R}$ with $n=8$, indicating that
the condition \eqref{eq:mu-cond} fails.

\subsection{Numerical tests with hexagonal lattice in Figure \ref{lattice3}}
\begin{figure}[hbt!]
\centering
\subfigure[$\mu=1$]{\label{loglog5_mu1_2 hex}
\includegraphics[width=8cm,trim={8cm 0.0cm 8cm 0.8cm},clip]{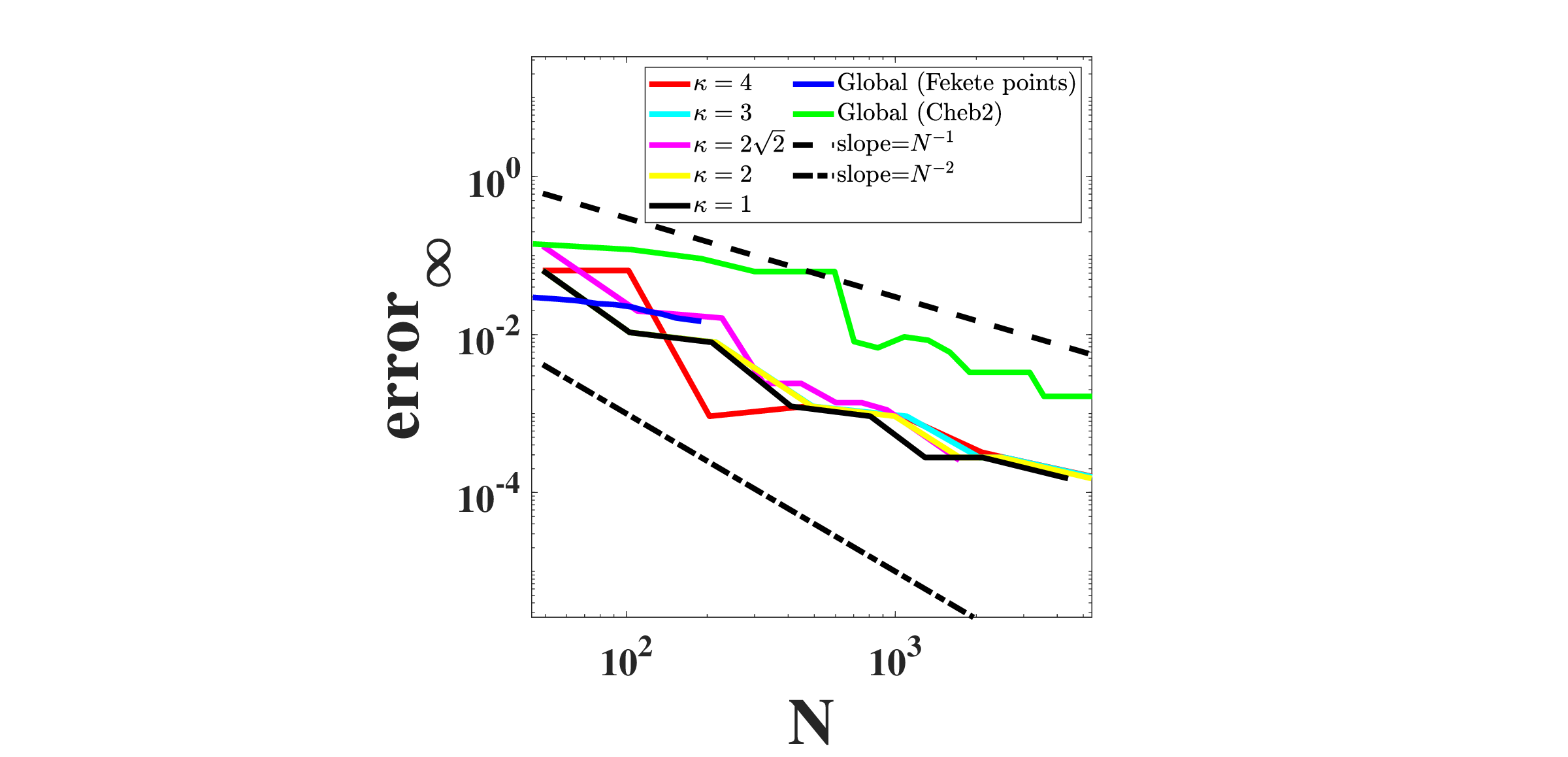}
}%
\quad
\subfigure[$\mu=0.5$]{\label{loglog5_mu2_2 hex}
\includegraphics[width=8cm,trim={8cm 0.0cm 8cm 0.8cm},clip]{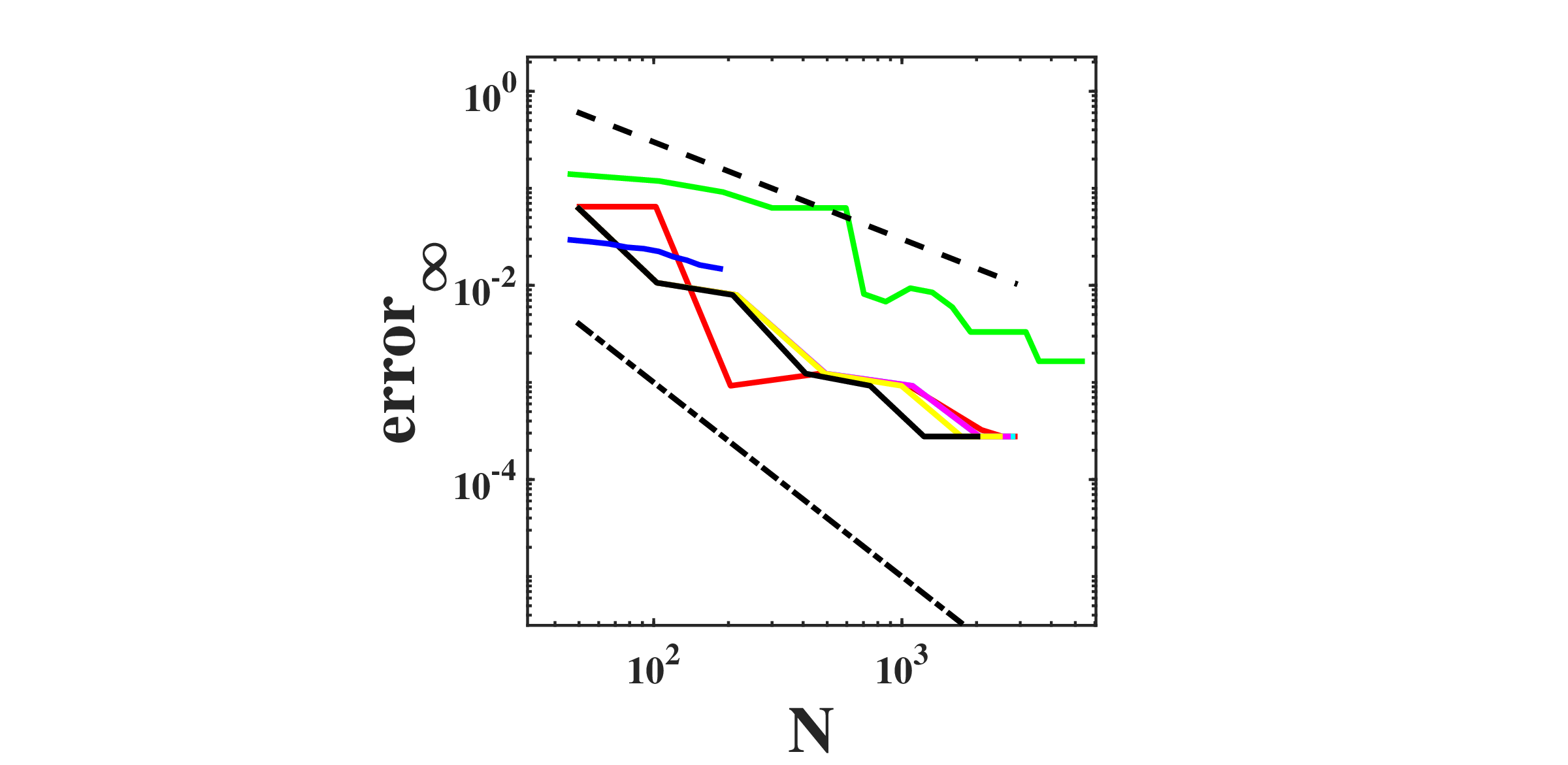}
}%
\quad
\centering
\caption{Maximum relative error with respect to the number of sampling points.}\label{loglog hex 1}
\end{figure}

\begin{figure}[hbt!]
\centering
\subfigure[loop 2]{\label{V5loop2 hex}
\includegraphics[width=3.7cm,trim={2cm 0.9cm 2cm 1.1cm},clip]{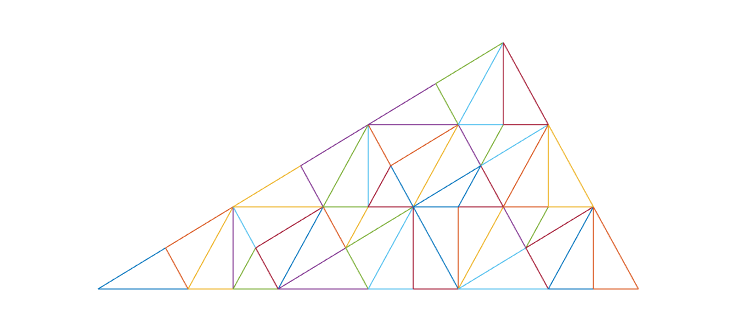}
}%
\quad
\subfigure[loop 4]{\label{V5loop4 hex}
\includegraphics[width=3.7cm,trim={2cm 0.9cm 2cm 1.1cm},clip]{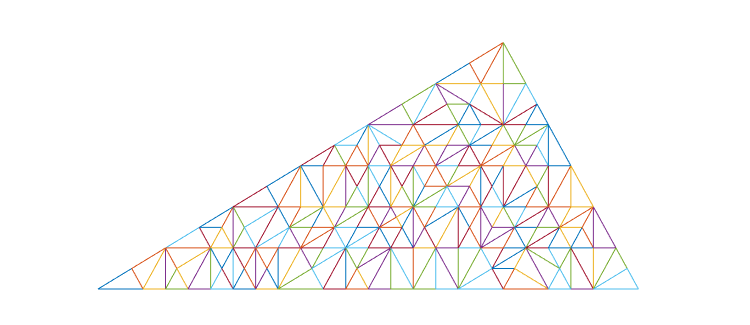}
}%
\quad
\subfigure[loop 6]{\label{V5loop6 hex}
\includegraphics[width=3.7cm,trim={2cm 0.9cm 2cm 1.1cm},clip]{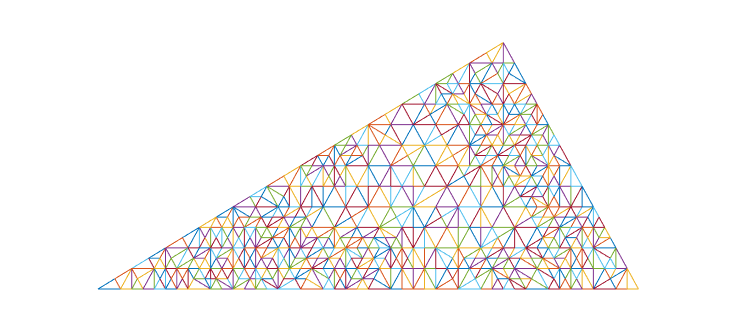}
}%
\quad
\subfigure[loop 8]{\label{V5loop8 hex}
\includegraphics[width=3.7cm,trim={2cm 0.9cm 2cm 1.1cm},clip]{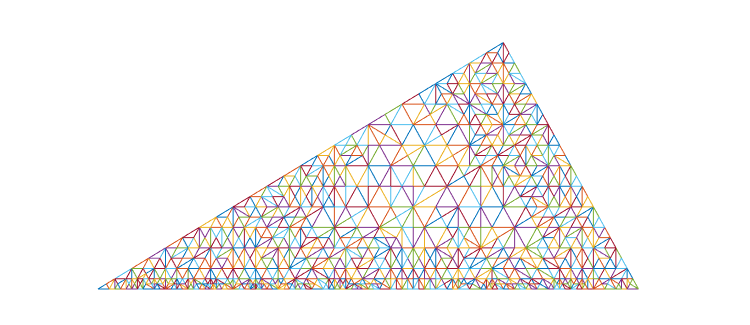}
}%
\quad
\subfigure[loop 2]{\label{V2loop2 hex}
\includegraphics[width=3.7cm,trim={2cm 0.9cm 2cm 1.1cm},clip]{figs_NEW/V5loop2_mesh.png}
}%
\quad
\subfigure[loop 4]{\label{V2loop4 hex}
\includegraphics[width=3.7cm,trim={2cm 0.9cm 2cm 1.1cm},clip]{figs_NEW/V5loop4_mesh.png}
}%
\quad
\subfigure[loop 6]{\label{V2loop6 hex}
\includegraphics[width=3.7cm,trim={2cm 0.9cm 2cm 1.1cm},clip]{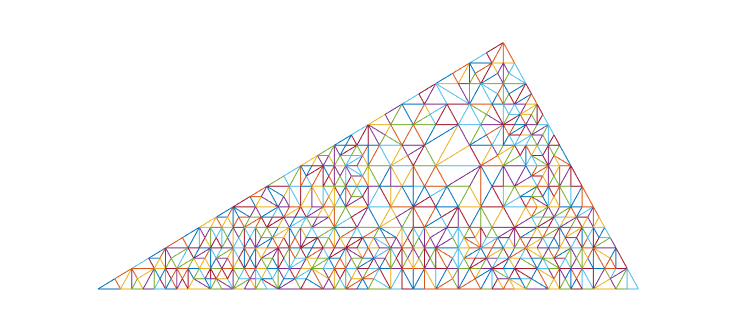}
}%
\quad
\subfigure[loop 8]{\label{V2loop8 hex}
\includegraphics[width=3.7cm,trim={2cm 0.9cm 2cm 1.1cm},clip]{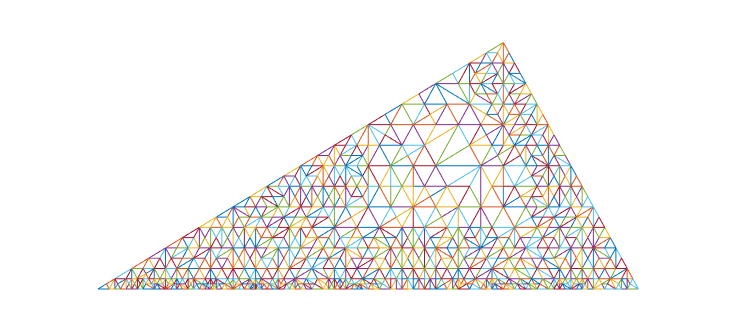}
}
\centering
\caption{Algorithm \ref{alg1} with $\kappa=2\sqrt{2}$ (the top row) and $\kappa=2$ (the bottom row).}\label{VXloopX hex}
\end{figure}

\begin{figure}[hbt!]
\centering
\subfigure[loop 8]{\label{hex kappa 2q2 error 1}
\includegraphics[width=8cm,trim={3cm 0.8cm 2cm 0.0cm},clip]{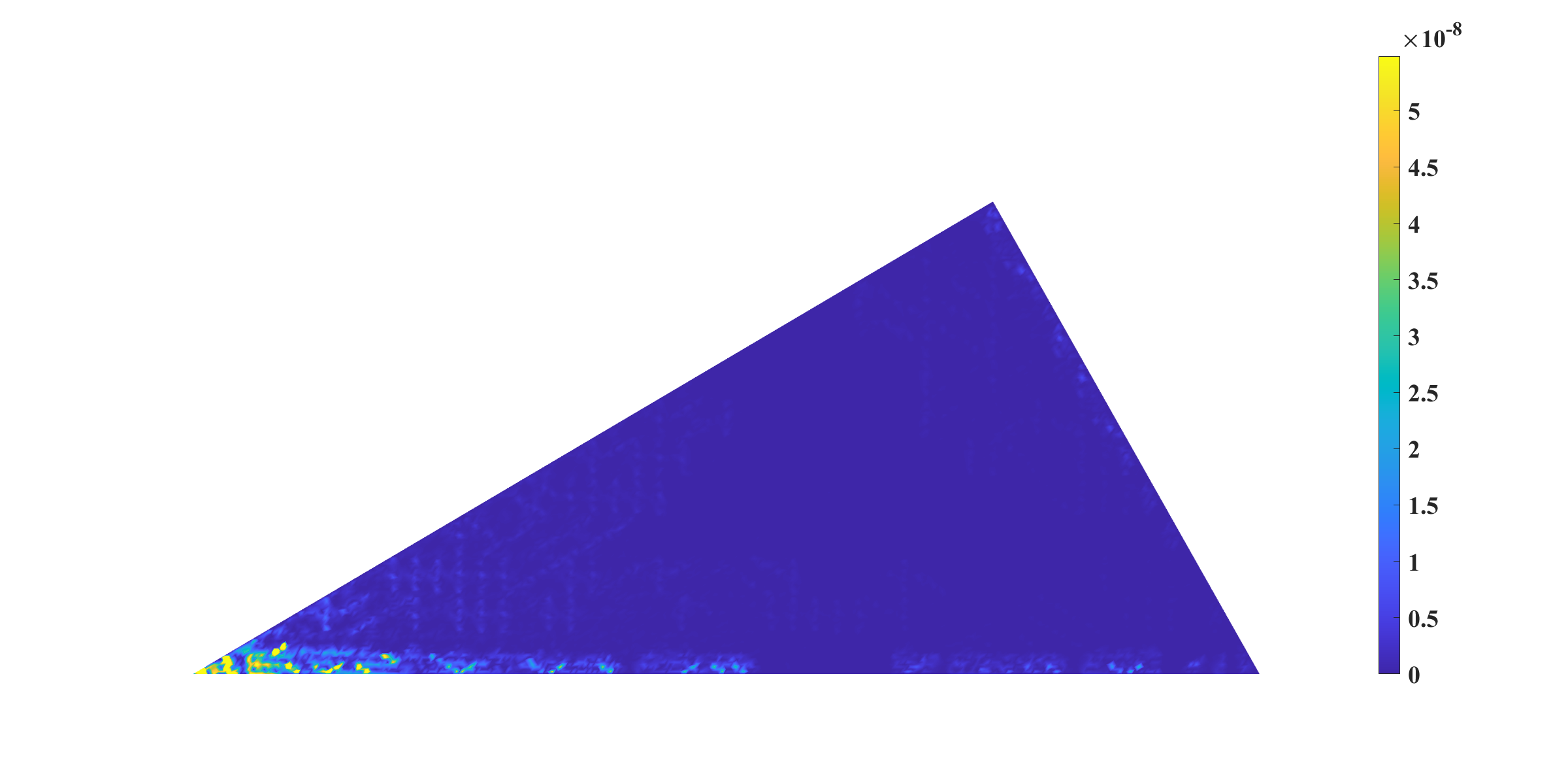}
}%
\quad
\subfigure[loop 8]{\label{hex kappa 2q2 error 2}
\includegraphics[width=8cm,trim={3cm 0.8cm 2cm 0.0cm},clip]{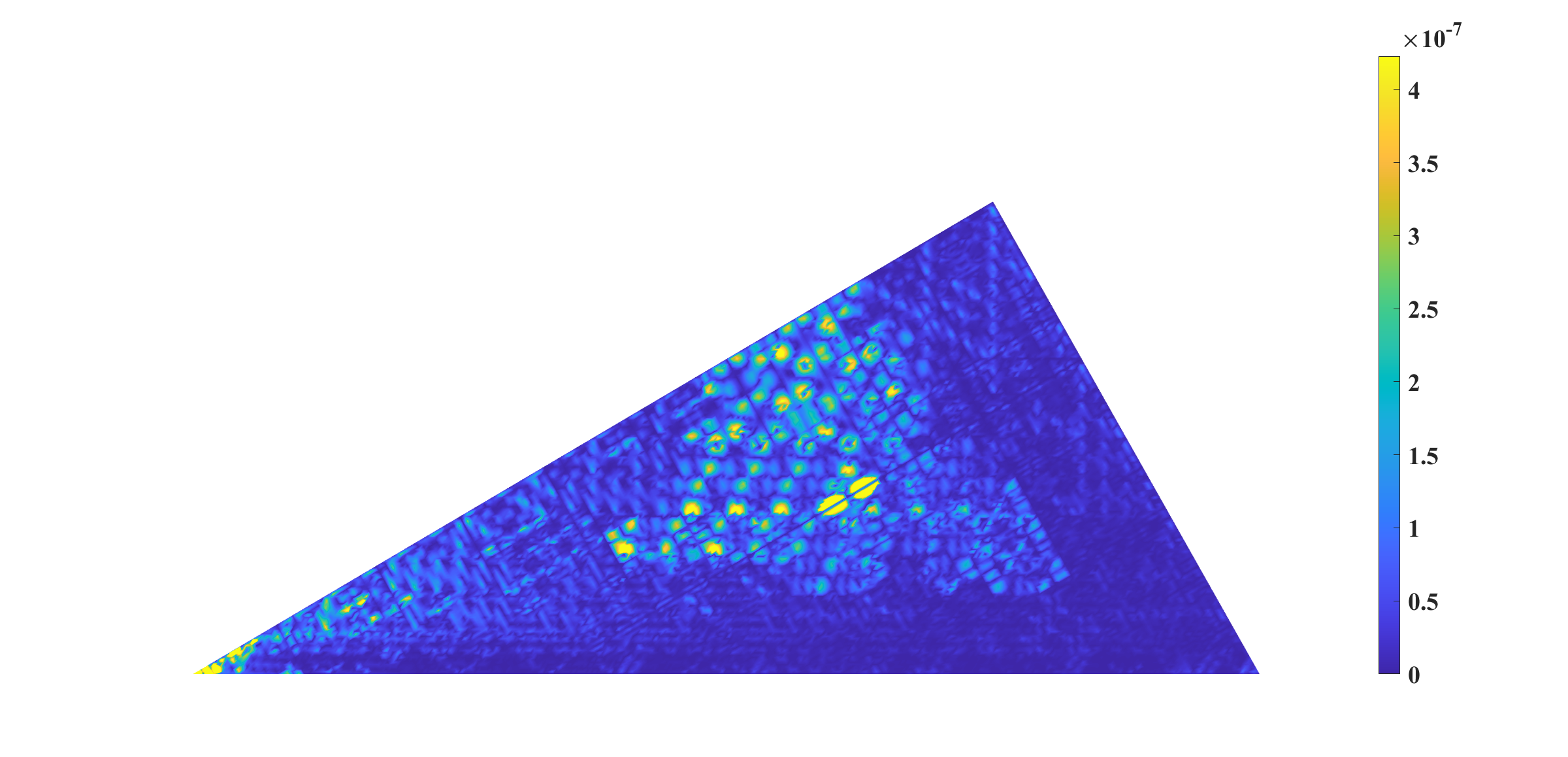}
}%
\caption{$\kappa=2\sqrt{2}$: pointwise relative error $e_i(\mathbf{k})$ under $\mu=1$ (left), $\mu=0.5$ (right).}\label{hex kappa 2q2 error}
\end{figure}

\begin{figure}[hbt!]
\centering
\subfigure[loop 8]{\label{hex kappa 2 error 1}
\includegraphics[width=8cm,trim={3cm 0.8cm 2cm 0.0cm},clip]{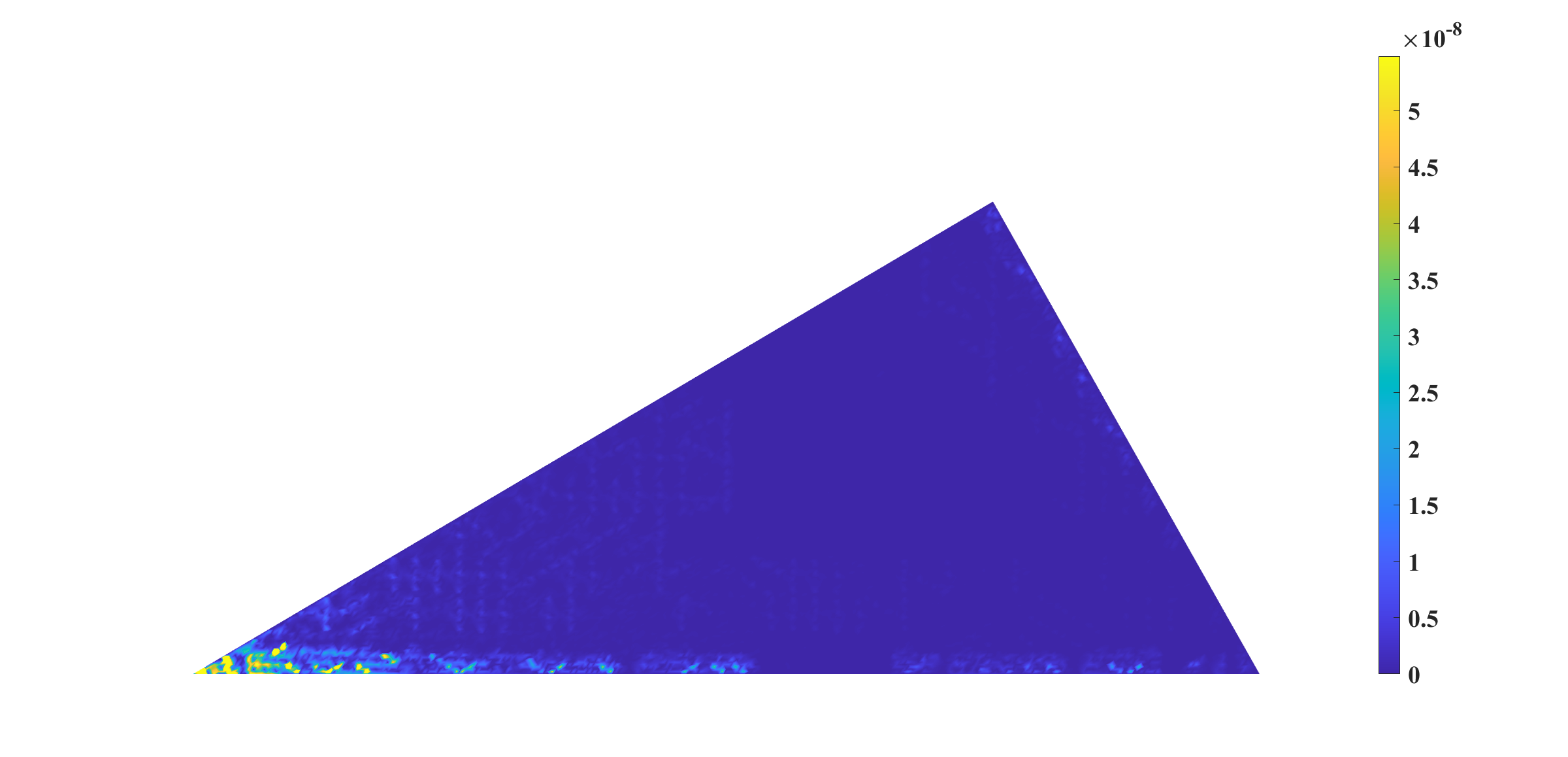}
}%
\quad
\subfigure[loop 8]{\label{hex kappa 2 error 2}
\includegraphics[width=8cm,trim={3cm 0.8cm 2cm 0.0cm},clip]{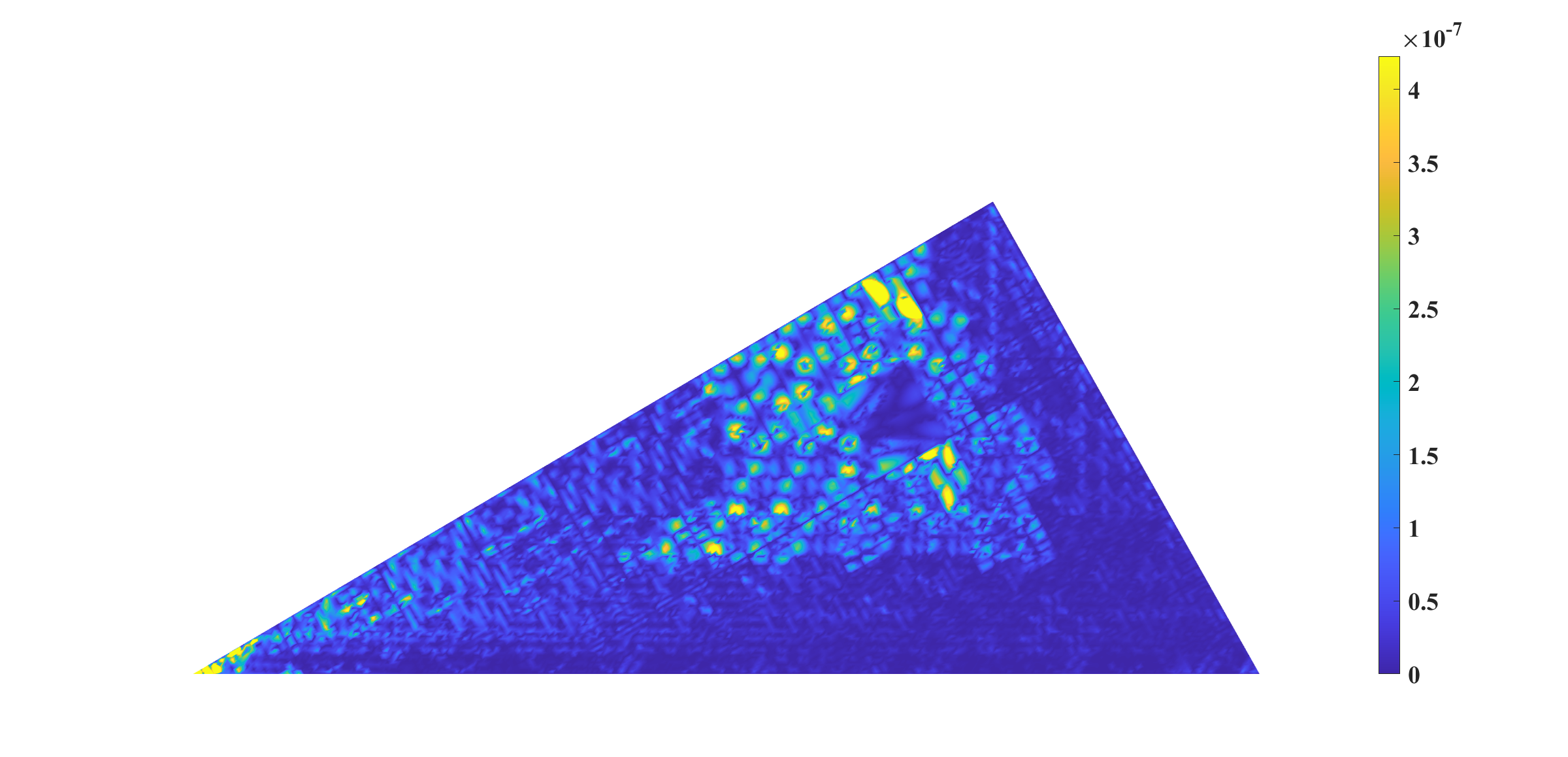}
}%
\caption{$\kappa=2$: pointwise relative error $e_i(\mathbf{k})$ under $\mu=1$ (left), $\mu=0.5$ (right).}\label{hex kappa 2 error}
\end{figure}

Now we study 2D PhCs with infinite periodic hexagonal lattice, cf. Figure \ref{lattice3}. The unit cell of the PhC is composed of six cylinders of dielectric material with dielectric constant $\epsilon=8.9$ embedded in the air and the lattice vectors $\mathbf{a}_1 = \frac{a}{4}({e}_1+\sqrt{3}{e}_2)$ and $\mathbf{a}_2 = \frac{a}{4}({e}_1-\sqrt{3}{e}_2)$ with lattice constant $a=3R$, $R$ is the length of hexagon edges and $r=\frac{1}{3}R$ is the radius of cylinders. We observe similar results as before. Figure \ref{loglog hex 1} shows the log-log plot of $\operatorname{error}_{\infty}$ against $N$. The adaptive mesh generated by Algorithm \ref{alg1} and pointwise relative errors $\operatorname{e}_i(\mathbf{k})$ with $\kappa=2,2\sqrt{2}$, $\mu=0.5,1$ are displayed in Figures \ref{VXloopX hex}, \ref{hex kappa 2q2 error}, \ref{hex kappa 2 error}, demonstrating the efficiency and reliability of Algorithm \ref{alg2}.

\section{Conclusion}\label{sec:conclusion}
In this work, we have presented an $hp$-adaptive sampling algorithm that can accurately approximate the band functions for two-dimensional photonic crystals. In the algorithm, we first generate an adaptive mesh that can refine elements containing singular points, and then assign proper polynomial degrees to each element based on its level. Additionally, we assign proper polynomial degrees to each edge to ensure a continuous global polynomial space. Finally, we develop an element-wise interpolation scheme using the global polynomial space to approximate band functions. We proved a first-order convergence rate in the worst-case scenario and presented several numerical tests to verify its performance. One interesting future research question is about approximating three-dimensional band functions, which is highly nontrivial due to a lack of provable local polynomial interpolation method on a tetrahedral domain. We plan to research on designing an optimal 3D PhCs with a maximum band gap.

\bibliographystyle{abbrv}
\bibliography{refer}
\end{document}